\documentclass[final,leqno,onefignum,onetabnum]{siamltex1213}
\usepackage{amsfonts,mathrsfs}
\usepackage{enumerate}

\newtheorem{remark}{{\em Remark}}[section]

\newtheorem{example}{{\em Example}}[section]

\def\3n{\negthinspace \negthinspace \negthinspace }
\def\2n{\negthinspace \negthinspace }
\def\1n{\negthinspace }

\def\dbE{\mathbb{E}}
\def\dbF{\mathbb{F}}

\def\dbP{\mathbb{P}}

\def\dbR{\mathbb{R}}

\def\={\buildrel \triangle \over =}

\def\ds{\displaystyle}

\def\ns{\noalign{\ss}}
\def\a{\alpha}
\def\b{\beta}

\def\e{\varepsilon}

\def\l{\lambda}
\def\m{\mu}
\def\n{\nu}
\def\si{\sigma}
\def\t{\tau}
\def\f{\varphi}
\def\th{\theta}
\def\o{\omega}

\def\G{\Gamma}
\def\D{\Delta}
\def\Th{\Theta}
\def\L{\Lambda}

\def\Om{\Omega}

\def\cA{{\cal A}}
\def\cB{{\cal B}}

\def\cD{{\cal D}}
\def\cE{{\cal E}}
\def\cF{{\cal F}}
\def\cG{{\cal G}}
\def\cH{{\cal H}}

\def\cK{{\cal K}}
\def\cL{{\cal L}}

\def\cU{{\cal U}}
\def\cV{{\cal V}}

\def\cl{{\cal l}}
\def\ss{\smallskip}

\def\mds{\medskip}

\def\q{\quad}
\def\qq{\qquad}
\def\hb{\hbox}

\def\lan{\mathop{\langle}}
\def\ran{\mathop{\rangle}}

\def\wt{\widetilde}

\def\cd{\cdot}

\def\ae{\hbox{\rm a.e.{ }}}
\def\as{\hbox{\rm a.s.{ }}}

\def\cl{\overline}

\def\({\Big (}
\def\){\Big )}
\def\[{\Big[}
\def\]{\Big]}

\def\mH{\mathscr{H}}

\def\bde{\begin{definition}}
\def\ede{\end{definition}}
\def\be{\begin{equation}}
\def\bel{\begin{equation}\label}
\def\ee{\end{equation}}
\def\bt{\begin{theorem}}
\def\et{\end{theorem}}
\def\bc{\begin{corollary}}
\def\ec{\end{corollary}}
\def\bl{\begin{lemma}}
\def\el{\end{lemma}}
\def\bp{\begin{proposition}}
\def\ep{\end{proposition}}
\def\bas{\begin{assumption}}
\def\eas{\end{assumption}}
\def\br{\begin{remark}}
\def\er{\end{remark}}
\def\ba{\begin{array}}
\def\ea{\end{array}}
\def\ed{\end{document}}

\def\eps{\epsilon}
\def\square#1{\vbox{\hrule\hbox{\vrule height#1%
     \kern#1\vrule}\hrule}}
\def\rectangle#1#2{\vbox{\hrule\hbox{\vrule height#1%
     \kern#2\vrule}\hrule}}

\font\tenbb=msbm10 \font\sevenbb=msbm7 \font\fivebb=msbm5

\newfam\bbfam
\scriptscriptfont\bbfam=\fivebb \textfont\bbfam=\tenbb
\scriptfont\bbfam=\sevenbb

\def\mr{\mathbb{R}}

\def\mrm{\mathbb{R}^m}
\def\mn{\mathbb{N}}

\def\mmf{\mathbb{F}}

\def\ms{\mathscr{S}}
\def\mb{\mathcal{B}}

\def\mx{\mathcal{X}}
\def\mt{\mathcal{T}}

\def\me{\mathbb{E}}
\def\mmu{\mathcal{U}_{ad}}
\def\bu{\bar{u}}

\def\T{T^{b}}
\def\eps{\varepsilon}
\def\t{\tau}

\newcommand{\inner}[2]{\left\langle#1,#2\right\rangle}

\title{Optimal control problems of forward-backward stochastic Volterra integral
equations with closed control regions}

\author{Tianxiao Wang\thanks{School of Mathematics, Sichuan University,
Chengdu, Sichuan Province, 610065, China. Email: xiaotian2008001@gmail.com.
This research of this
author was supported by NSF of China under grant
11231007, 11301298, 11401404 and 11471231, China Postdoctoral Science Foundation (2014M562321).}\and
Haisen Zhang\thanks{Corresponding author. School of Mathematics and Statistics,
Southwest University, 400715 Chongqing, China. Email: haisenzhang@yeah.net. The research of this
author is partially supported by Natural Science Foundation of China under grant 11471231, the fundamental research funds for the central universities under grants SWU114074 and XDJK2015C142, and the Natural Science Foundation Project of ChongQing  CSTC under grant 2015jcyjA00017.}}

\begin{document}
\maketitle
\slugger{sicon}{xxxx}{xx}{x}{x--x}

\begin{abstract}
Optimal control problems of forward-backward stochastic Volterra
integral equations (FBSVIEs, in short) with closed control
regions are formulated and studied. Instead of using spike variation method
as one may imagine, here we turn to treat the non-convexity of the control regions by borrowing some tools in set-valued analysis and adapting them into our stochastic control systems. A duality principle between linear backward stochastic Volterra integral equations and linear stochastic Fredholm-Volterra integral equations with conditional expectation are derived, which extends and improves the corresponding results in \cite{Shi-Wang-Yong 2015}, \cite{Yong 2008}. Some first order necessary optimality conditions for optimal controls of FBSVIEs are established. In contrast with existed common routines to treat the non-convexity of stochastic control problems, here only one adjoint system and one-order differentiability requirements of the coefficients are needed.
\end{abstract}



\pagestyle{myheadings}
\thispagestyle{plain}
\markboth{T.~Wang and H.~Zhang}{Optimal control problems of FBSVIEs}

\section{Introduction}

Let $T>0$ and $(\Omega,\cF, \mmf, \dbP)$  be a complete filtered
probability space (satisfying the usual conditions), on which a
$1$-dimensional standard Wiener process $W(\cdot)$ is defined such
that $\mmf=\{\cF_{t} \}_{0\le t\le T}$ is the natural filtration
generated by $W(\cdot)$ (augmented by all of the $\dbP$-null sets).
Consider the following controlled stochastic differential equation
\begin{equation}\label{State-SDEs}
\left\{
\begin{array}{l}
dX(t)=b(t,X(t),u(t))dt+\sigma(t,X(t),u(t))dW(t),\ \ \ t\in[0,T],\\
X(0)=x_0,
\end{array}\right.
\end{equation}
with cost functional
\begin{equation}\label{Cost-SDEs}
J(u(\cdot))=\me\Big[\int_{0}^{T}f(t,X(t),u(t))dt+h(X(T))\Big].
\end{equation}
Here $u(\cdot)$ is the control variable valued in the control region $U\subset
\dbR^{l}$, $X(\cdot)$ is the state variable valued in $\mr^n$ (for
some $n\in \mn$), and $b,\sigma$, $f$ and $h$ are given functions.
The stochastic optimal control problem is to find a control variable
$\bar{u}(\cdot)$ belonging to the admissible control set $\mmu$ (which will be defined later) such that
$$J(\bar{u}(\cdot))=\inf_{u(\cdot)\in\mmu}J(u(\cdot)).$$

For above problem, one of the central topics is to establish the
necessary conditions for optimal controls. Many contributions in
this field were made ever since the work of \cite{Kushner-Schweppe},
see e.g. \cite{Bensoussan 1981}, \cite{Bismut 1978}, \cite{Haussmann
1976} and references cited therein. However, the general case with
control-dependent diffusion term and non-convex control region were
untouched until the work \cite{Peng90}. Besides the
standard spike variation method and the useful tool of It\^{o}
formula, another indispensable notion in \cite{Peng90} is the introduced
second-order adjoint equation. Note that the later one is actually a
linear backward stochastic differential equation (BSDE, in short),
the point of which also reveals the crucial role of BSDEs in
stochastic optimal control problems. In fact, besides their
wide use in the stochastic control problems, BSDEs are also
applicable in other areas, such as mathematical finance. For example, a large class of risk
measures or stochastic differential utility can be represented by
the solutions of proper BSDEs (see \cite{Duffie-Epstein
1992}, \cite{Gianin 2006}). Moreover, according to e.g. \cite{Peng-Wu 1999}, \cite{Wang-Wu-Xiong 2013}, \cite{Yong 2010}, the financial/economic applications largely motivate people to study the optimal control problem for forward-backward stochastic differential equations (FBSDEs, in
short),
\bel{State-FBSDEs}~~~~~\left\{\ba{ll}
\ns\ds X(t)\!=\!x_0\!+\!\int_0^t\!b(s,X(s),u(s))ds\!+\!\int_0^t\!\sigma(s,X(s),u(s))dW(s),\
t\!\in\![0,T],\\
\ns\ds Y(t)\!=\!h(X(T))\!+\!\!\int_t^T\!\!\!g(s,X(s),Y(s),Z(s),u(s))ds\!-\!\!\int_t^T\!\!\!Z(s)dW(s),\  t\!\in\![0,T],
\ea\right. \ee
associated with cost functional
\bel{Cost-FBSDEs}\ba{ll}
\ns\ds
J(u(\cdot))=\me\Big[\int_{0}^{T}f(t,X(t),u(t),Y(t),Z(t))dt+h(x(T),Y(0))\Big].
\ea\ee
However, such control problem with general control region also kept still for nearly decade until
recently \cite{Wu 2013} and \cite{Yong 2010} gave some excellent solutions along this.
Actually, the reliance of diffusion term on control variable and the
limited integrability of process $Z(\cd)$ makes the second order
Taylor-type expansion becomes impossible, not to mention the
deriving of maximum principle. In order to get around these
essential difficulties, the authors in \cite{Wu 2013} and \cite{Yong
2010} transformed equivalently the original forward-backward problem
into a new forward control system case with initial-terminal state constraints. Then
the original issue can be solved by working on the later one.

In this paper we study the
optimal control problem of forward-backward stochastic
Volterra integral equations (FBSVIEs, in short),
\bel{State-FBSVIEs}~~\left\{\2n\1n\ba{ll}
\ns\ds
X(t)\!=\!\f(t)\!+\!\int_0^t\!b(t,s,X(s),u(s))ds\!+\!\int_0^t\!\si(t,s,X(s),u(s))dW(s), \\
\ns\ds
Y(t)\!=\!\psi(t,X(T))\!+\!\int_t^T\!\!g(t,s,X(s),Y(s),Z(t,s),u(s))ds
\!-\!\int_t^T\!\!Z(t,s)dW(s).\\
\ea\right.\ee
Mathematical speaking, FBSVIE (\ref{State-FBSVIEs}) is an extension of FBSDE (\ref{State-FBSDEs}). The motivations of
our study are based on the following aspects. To begin with, let us take more closer glances at both equations
(\ref{State-SDEs}) and (\ref{State-FBSDEs}), from which one can
understand a fundamental structure of differential systems: time
consistency (or semi-group property). Actually it is just this
inherent feature that makes some mathematical treatments, such as
the well-known dynamic programming principle (\cite{Yong-Zhou
1999}), or the dynamic risk measures by BSDEs (\cite{Gianin 2006}),
applicable and useful. However, from practical point of view, such
character seems to make the described system rather ideal,
even in deterministic setting.  For example, the physical
meaningfulness of heat equation has been doubted due to its property
of infinite speed of propagation. To solve this problem, one common
way is to add some memory effects into this partial differential
equations, see e.g. \cite{Yong-Zhang 2011}. Inspired by this point,
we would like to replace the forward equation of (\ref{State-FBSDEs}) with some
stochastic system with memory, like stochastic Volterra integral equation
(SVIE, in short). For the risk measures/differential utility
represented by BSDEs, some recent study (e.g.
\cite{Ekeland-Mbodji-Pirvu 2012}, \cite{Wang-Yong 2015}) also
indicates a tendency to replace them with general time inconsistent
counterparts, like the ones via backward stochastic Volterra integral equations
(BSVIEs, in short). Therefore, here we use the controlled BSVIEs in (\ref{State-FBSVIEs}) instead of
the controlled BSDE in (\ref{State-FBSDEs}). Furthermore arguments along this can be found in \cite{Shi-Wang-Yong 2015}.
Here we would like to mention some related study on SVIEs and
BSVIEs, e.g. \cite{Agram-Oksendal 2015}, \cite{Bonaccorsi-Confortola-Mastrogiacomo 2012}, \cite{Lin 2002}, \cite{Oksendal-Zhang 2010}, \cite{Pardoux-Protter 1990}, \cite{Protter 1985}, \cite{Shi-Wang-Yong
2015}, \cite{Wang-Yong 2015}, \cite{Wang-Zhang 2007}, \cite{Yong 2008}.

The purpose of this paper is to establish some first order necessary optimality conditions for FBSVIEs with non-convex control region.
If one follows the conventional approaches in FBSDEs case (e.g. \cite{Peng-Wu 1999}, \cite{Wang-Wu-Xiong 2013}, \cite{Wu 2013}, \cite{Yong 2010}), there
are two essential difficulties one has to face with. In the first place, the transformation between controlled FBSDEs and another controlled SDEs with state constraints appearing in \cite{Wu 2013} and \cite{Yong 2010} actually made use of time consistency of differential systems, which of course does not fit in the SVIEs framework. Second, the inherent structure of Volterra integral systems makes the duality between SDEs and BSDEs via It\^{o}'s formulation no longer work well for SVIEs and FBSVIEs (see \cite{Shi-Wang-Yong 2015}), hence many developed excellent tricks in FBSDEs case are absent here. Therefore, we need to provide more efficient techniques rather than following the traditional ones.

In contrast with spike variation, in this paper we will use a quite different variational technique to deal with the non-convexity of the control region $U\subset\dbR^l$. To show the basic ideas involved, let us firstly recall the convex case. Let $\bar u(\cd)$ be  an optimal control and define
$$
\cV:=\{v(\cd)=[u(\cd)-\bar u(\cd)]\in L^2_{\dbF}(0,T;\dbR^l)\ |\  u(\cd)\in\cU_{ad}\}.
$$
When $U$ is convex, it is clear that, for any $v(\cdot)\in \cV$ and $\e\in(0,1)$,
\bel{Perturbation-manner-1}\ba{ll}
\ns\ds
\bar u(\cd)+\e v(\cd)=\Big[\bar u(\cd)+\e\big(u(\cd)-\bar u(\cd)\big)\Big]\in\cU_{ad}.
\ea\ee
Such kind of perturbation is named a convex variation of $\bar u(\cd)$. Based on this result, by introducing suitable variational equation, adjoint equation and related duality skills one can obtain the required necessary optimality condition. Here $\cV$ can be seen as the set of perturbation direction for $\bar u(\cd)$. However, when the control region $U$ is non-convex, there may exists $u_1(\cd)$ such that for any $\e>0$, $\bar u(\cd)+\e\big(u_1(\cd)-\bar u(\cd)\big)$ does not belong to $\cU_{ad}$, i.e., $\cV$ is a little bit of large as a set of perturbation direction in this case. So we should find another suitable way to choose the set of perturbation direction but without losing the basic procedures from variational equation to duality principle. One way to do so is to find set $\bar\cV\subset L^2_{\dbF}(0,T;\dbR^l)$ such that for any $v(\cd)\in\bar\cV$ and small $\e>0$, there exists $v_{\e}(\cd)\in L^2_{\dbF}(0,T;\dbR^l)$ satisfying
\bel{Perturbation-manner-2}\ba{ll}
\ns\ds
\bar u(\cd)+\e v_{\e}(\cd)\in\cU_{ad},\qq \dbE\int_0^T|v_\e(s)-v(s)|^2ds\rightarrow0,\ \ \e\rightarrow0.
\ea\ee
In some sense condition (\ref{Perturbation-manner-2}) can be seen as a extension/relaxtion of (\ref{Perturbation-manner-1}).
However, the question is: do the set $\bar\cV$ and following-up tricks exist? Fortunately, as we will see next, the adjacent cone of $\cU_{ad}$ at $\bar u(\cd)$ (see Definition \ref{Definition 2.2}) is a good choice of $\bar\cV$ satisfying (\ref{Perturbation-manner-2}). Actually, such perturbation approach via the adjacent cone is called variational analysis approach in the literature. Note that the variational analysis approach has a long history and been used extensively in optimization and optimal control theory under the deterministic setting, see the book \cite{Bonnans00}. Using this method, \cite{Hoehener12} established a second order integral type necessary condition for optimal control problem of ordinary differential equations with state constraint, which was later improved into the pointwise form in \cite{Frankowska13} with delicate analysis. For the stochastic case, \cite{Bonnans12} firstly used the variational analysis approach to deal with controlled SDEs and obtained both the first and second order
integral type necessary condition with convex and closed control regions. In this paper, we adopt such a variational analysis approach under
the FBSVIEs setting with closed control region (but not necessary to be convex) and obtain some pointwise necessary conditions of optimal controls. Notice that the pointwise form seems more appropriate than the integral counterparts in \cite{Bonnans12} from the view of mathematical control theory.

At this moment we would like to point out some novelties of following-up studies.  In the first place, when the control region is nonconvex, compared with existed spike variation method it seems that our variational approach (under some structural assumptions on $U$, see Remark \ref{strucural assum}) is more advantageous. Actually, in the particular SDEs case, to derive the maximum principle the coefficients in the traditional literature are required to be twice differentiable with respect to variable $x$. Moreover, besides the first-order adjoint equations, the second-order adjoint equations are also indispensable in the stochastic setting (see \cite{Peng90}). Nevertheless, here we only need the differentiability of the coefficient up to the first order. In addition, only one adjoint equation is needed even though the control region is allowed to be non-convex. In the second place, since we are using set-valued analysis in the stochastic framework, it seems like the procedures of transforming the integral necessary condition into the pointwise form are essentially different from the existed counterparts. In fact, similar transformation can be directly derived via contradiction arguments if we use spike/convex variational method. However, some new features arise under our framework which make the pointwise process by no means straightforward as before. Fortunately, by borrowing some existed well-known results in the set-valued analysis (see \cite{Aubin90}) the pointwiseness arguments can be successfully done. To our best such result, i.e. Lemma \ref{Lemma-integral-pointwise} appears for the first time. Furthermore, its generality also makes it applicable in obtaining necessary optimality conditions for some other stochastic control systems. In the third place, when we are trying to use the existed results on set-valued analysis (see \cite{Aubin90}), the incompleteness of product measure space $(\Omega\times[0,T],\cF_{T}\otimes\cB([0,T]),\l\times\dbP)$ under stochastic setting does not fit their completeness requirement (see e.g. Proposition \ref{Pro-adjacent} next). Hence further works also need to be done to fill the gap between the two. Since such a problem can be avoided for the study on deterministic system, like \cite{Frankowska13}, \cite{Lou10}, it thus implies a new distinction between stochastic control problem and deterministic case.

The rest of this paper is organized as follows. In Section 2, we
list some notations, spaces and preliminary results. In Section 3,
we introduce the main results of this paper and give some examples.
Finally, in Section 4 we give the proof of our main result, as well as a general dual principle and some pointwise procedures.

\section{Some preliminaries}

In the first place let us introduce some notations. For
$H=\dbR^n,\dbR^{n\times m}$, etc., we denote its norm by $|\cd|$.
For $0\le s<t\le T$, we define
$$\ba{ll}
\ns\ds L^2_{\cF_t}(\Om;H):=\Big\{\xi:\Om\to H\bigm|\xi\hb{ is
$\cF_t$-measurable, }\dbE|\xi|^2<\infty\Big\},
\ea$$
$$\ba{ll}
\ns\ds L^2_{\cF_T}(s,t;H):=\Big\{X:[s,t]\times\Om\to H\bigm|X(\cd)\hb{
is $\cF_T\otimes\cB([s,t])$-measurable},\\
\ns\ds\qq\qq\qq\qq\qq\qq\qq\qq\qq\qq\q
\hb{such that}\ ~\dbE\int_s^t|X(r)|^2dr<\infty\Big\},\\
\ns\ds L^2_{\cF_T}\big(\Om;C([s,t];H)\big):=\Big\{X:[s,t]\times\Om\to
H\bigm|X(\cd)\hb{ is $\cF_T\otimes\cB([s,t])$-measurable, }\\
\ns\ds\qq\qq\qq\qq\qq\qq\q\hb{and has  continuous paths, } \dbE\(\sup_{r\in[s,t]}|X(r)|^2\)<\infty\Big\},\\
\ns\ds C_{\cF_T}\big([s,t];L^2(\Om;H)\big)\1n:=\1n\Big\{X\1n:\1n[s,t]\to
L^2_{\cF_T}(\Om;H)\bigm|X(\cd)\hb{ is
continuous from $[s,t]$ to }\\
\ns\ds\qq\qq\qq\qq\qq\qq\qq\qq\qq\qq L^2_{\cF_T}(\Omega,H),\ \ \sup_{r\in[s,t]}\dbE|X(r)|^2<\infty\Big\}.\ea$$
Also, we define
$$L^2_\dbF(s,t;H):=\Big\{X(\cd)\in L^2_{\cF_T}(s,t;H)\bigm|X(\cd)\hb{ is $\dbF$-adapted}\Big\}.$$
The spaces $L^2_\dbF\big(\Om;C([s,t];H)\big)$ and
$C_\dbF\big([s,t];L^2(\Om;H)\big)$ can be defined in the similar way.
Further, we denote
$$\ba{ll}
\ns\ds\D:=\Big\{(t,s)\in[0,T]^2\Bigm|t\le
s\Big\},\q\D^*:=\Big\{(t,s)\in[0,T]^2\Bigm|t\ge s\Big\}\equiv
\cl{\D^c},\\
\ns\ds\cH^2_\D(0,T;H):= L^2_\dbF(0,T;H)\times
L^2_\dbF(\D;H),\\
\ns\ds \cH^2(0,T;H):= C_{\dbF}([0,T];L^2(\Omega;H))\times
C\big([0,T];L^2_\dbF(0,T;H)\big),\\
\ns\ds \cH^2_{1}(0,T;H):= L^2_{\dbF}(0,T;H)\times
L^2\big(0,T;L^2_\dbF(0,T;H)\big),
\ea
$$
where
$$\ba{ll}
\ns\ds L^2_\dbF(\D;H):=\Big\{Z:\D\times\Om\to H\bigm|s\mapsto
Z(t,s)\hb{ is $\dbF$-adapted on $[t,T]$, $\ae t\in[0,T]$,}\\
\ns\ds\qq\qq\qq\qq\qq\qq\qq \|Z(\cd\,,\cd)\|_{L^2_\dbF(\D;H)}^2:=\dbE\int_0^T
\int_t^T|Z(t,s)|^2dsdt<\infty\Big\},\\
\ns\ds C([0,T];L^2_\dbF(0,T;H)):=\Big\{Z:[0,T]^2\times\Om\to
H\bigm|s\mapsto Z(t,s)\hb{ is $\dbF$-adapted}\hb{ on $[0,T]$,}\\
\ns\ds\qq\qq\qq\qq\qq\q \hb{ $\forall t\in[0,T]$,}\ t\mapsto Z(t,\cd)\hb{ is continuous in}\ L^2_\dbF(0,T;H)\ \hb{and}\\
\ns\ds\qq\qq\qq\qq\qq\ \ \|Z(\cd\,,\cd)\|_{C(0,T;L^2_\dbF(0,T;H))}^2:=
\sup_{t\in[0,T]}\dbE\int_0^T|Z(t,s)|^2ds<\infty\Big\},\\
\ns\ds L^2(0,T;L^2_\dbF(0,T;H)):=\Big\{Z:[0,T]^2\times\Om\to
H\bigm|s\mapsto Z(t,s)\hb{ is $\dbF$-adapted}\hb{ on $[0,T]$,}\\
\ns\ds\qq\q\hb{ $t\in[0,T]$,\ a.e. \ }\|Z(\cd\,,\cd)\|_{L^2(0,T;L^2_\dbF(0,T;H))}^2:=
 \dbE\int_0^T\int_0^T|Z(t,s)|^2dsdt<\infty\Big\}.\ea$$
For readers' convenience, next let us introduce some necessary notations
and results on set-valued analysis, see \cite{Aubin90}. Let $X$ be a
Banach space with norm $\|\cdot \|$. For any $x\in X$, denote $B(x,\eps):=\{y\in X\ |\ \|y-x\|< \eps\}$ with its closure $\bar{B}(x,\eps)$. For any subset $K\subset X$, denote by $int K$, $cl K$ and $bd K$ the interior, closure and boundary of $K$. The
distance between a point $x$ and a closed set $K$ in $X$ is defined by
$\hb{dist}(x,K):= \inf_{y\in K} \|y-x\|$. The following definition can be
found in Definition 4.1.5 (p. 126) of \cite{Aubin90}.

\medskip
\bde \label{Definition 2.2} \rm Let $K$ be a closed subset of a
Banach space $X$, $x\in K$. The adjacent cone $\T_{K}(x)$ is defined
by
$$\T_{K}(x):=\Big\{v\in X\ \Big|\ \lim_{h\to 0^+} \frac{\hb{dist}(x+hv,K)}{h}=0\Big\}.$$
\ede
Obviously, for any $x\in int K$, $\T_{K}(x)=X$, and for any $x\in bd K$, $\T_{K}(x)$ is a closed set and  $0\in \T_{K}(x)$. When $K$ is convex,
$$\T_{K}(x)=cl\Big\{\alpha(y-x)\ \Big|\ \alpha\ge0,\ y\in K \Big\}.$$

Now let us introduce some characterizations of
the adjacent cone which can be found in page 128 of \cite{Aubin90}.

\medskip
\bl\label{equivdef for adjacent cone}  \rm Let $K$ be a closed subset of a
Banach space $X$, $x\in K$. The following assertions are equivalent.
\begin{enumerate}[\rm(i)]
\item $v\in \T_{K}(x)$.
\item  For any $\eps>0$, there exists an $\alpha>0$ such that for any $h\in (0,\alpha)$ one can find a vector $v_{h}\in B(v,\eps)$ so that  $x+hv_{h}\in K$.
\item  For any $h_{n}\to 0^+$, there exists a sequence $\{v_{n}\}_{n=1}^{\infty}$ such that $v_{n}\to v$ in $X$ as $n\to \infty$  and for any $n$, $x+h_{n}v_{n}\in K$.
\end{enumerate}
\el

\rm
\medskip
\begin{remark}\label{remark2-1}
\em
Note that for any fixed $h>0$, there exists $y_{h}\in K$ such
that
$$\|y_{h}-x-hv \|\le \inf_{y\in K} \|y-x-hv\|+h^2.$$
If $v\in \T_{K}(x)$, by denoting $v_{h}:=\frac{y_{h}-x}{h}$, it
follows that $v_{h}\to v$ as $h\to 0^+$, and
$x+hv_{h}\in K$. Such a point will be useful in the sequel.
\end{remark}
\medskip

In the following, let $(\Xi, \ms, \mu)$ be a
$\sigma$-finite measure space, $X$ a separable Banach space,
$F:\Xi\rightsquigarrow X$ a set-valued map. For any $\xi\in \Xi$,
$F(\xi)$ is called the value (or the image) of $F$ at $\xi$. The
domain of $F$ is the subset of $\xi\in \Xi$ such that $F(\xi)$ is
not empty, i.e.,
$$\hb{Dom}(F):=\{\xi\in \Xi\ |\ F(\xi)\neq \emptyset\}.$$
The image of $F$ is defined by
$$\hb{Im}(F):=\bigcup_{\xi\in \Xi}F(\xi). $$
For map $F$, let us define its graph $\hb{Graph}(F)$ the subset of
product space $\Xi\times X$ as,
$$ \hb{Graph}(F):=\big\{(x,y)\in\Xi\times X\big| y\in F(x)\big\}.
$$
Suppose $F:\Xi\rightsquigarrow X$ is a set-valued map with closed
image. $F$ is called measurable if for any Borel set $A \subseteq X$,
the inverse image of $F$ is measurable, i.e.,
$$F^{-1}(A):=\{\xi\in \Xi\ |\ F(\xi)\cap A\neq \emptyset\}\in \ms.$$
Note that the domain of a measurable map is measurable as well as
its complement $\big\{\xi\in\Omega\big| F(\xi)=\emptyset\big\}.$ The
following result gives one criteria for the measurability, i.e.
Theorem 8.1.4 (p.310) of \cite{Aubin90}.

\medskip
\bp \label{Pro-measurable}
 \rm Let $(\Xi,\ms,\mu)$ be a complete $\sigma$-finite measure space, $X$
a complete separable metric space and $F$ a set-valued map from
$\Xi$ to $X$ with nonempty closed images. Then $F$ is measurable if
and only if the graph of $F$ belongs to $\ms\otimes\mx$, where $\mx$
is the Borel $\sigma$-algebra of $X$.
\ep

\medskip
\rm
In what follows, we also need the notion of
measurable selection of a given set-valued map, see Definition 8.1.2
(p.308) of \cite{Aubin90},

\medskip
\bde \label{Def-mea-sel}\rm
Let $(\Xi,\ms)$ be a measurable space and $X$ a complete separable
metric space. Consider a set-valued map $F$ from $\Xi$ to $X$. A
measurable map $f:\Xi\to X$ satisfying $f(\xi)\in F(\xi)$ for
any $\xi\in\Xi$, is called a measurable selection of $F$.
\ede

\medskip
The following result comes from Theorem 8.1.3 (p.308) in
\cite{Aubin90},

\medskip
\bp \label{Pro-mea-sel}
 \rm
Let $X$ be a complete sparable metric space, $(\Xi,\ms)$ a
measurable space, $F:\Xi\rightsquigarrow X$ a measurable set-valued map with nonempty
closed values. Then there exists a measurable
selection of $F$.
\ep
\rm

\medskip
Eventually, let us look at one result on the adjacent cone, see also
Theorem 8.5.1 (p.324) of \cite{Aubin90}.

\medskip
\bp \label{Pro-adjacent}
 \rm Suppose $(\Xi,\ms,\mu)$ is a complete $\sigma$-finite measure
space, and $X$ is a separable Banach space, $U\subseteq X$ is a
closed set. Then for any $k(\cd)\in\cK$, with
$$
\cK:= \big\{k(\cd)\in L^p(\Xi,\ms,\mu)\big| \ \hb{for almost
all} \ \xi\in\Xi,\ k(\xi)\in U\big\},
$$
and $p\geq1$, the set-valued map $\T_{U}(k(\cd))$ is
$\ms$-measurable, and $\mt_k\subseteq\T_{\cK}(k(\cd))$ where $\mt_k$ is defined as
$$\mt_k := \big\{l(\cd)\in L^p(\Xi,\ms,\mu)\big|\
\hb{for almost all} \ \xi\in\Xi,\ \ l(\xi)\in\T_{U}(k(\xi))\big\}.
$$
\ep

\rm

\section{Optimal Control Problems and Maximum Principles}

Let us recall the controlled froward-backward stochastic Volterra integral equations (FBSVIEs, in short):
\bel{FBSVIE-3-1}~~\left\{\2n\1n\ba{ll}
\ns\ds
X(t)\!=\!\f(t)\!+\!\int_0^t\!b(t,s,X(s),u(s))ds\!+\!\int_0^t\!\si(t,s,X(s),u(s))dW(s), \\
\ns\ds
Y(t)\!=\!\psi(t,X(T))\!+\!\int_t^T\!\!g(t,s,X(s),Y(s),Z(t,s),u(s))ds\!-\!\int_t^T\!\!Z(t,s)dW(s),
\ea\right.\ee
where $t\in[0,T]$, and $u(\cd)$ belongs to the set of admissible controls $\cU_{ad}$ defined by
$$\cU_{ad}:=\Big\{u(\cdot )\in L_{\dbF}^2(0,T;\dbR^\ell)\bigm|u(t)\in
U,~\ae t\in[0,T],~\as\2n\Big\},$$
with $U$ being a nonempty closed subset of $\dbR^\ell$.

\mds

\bde \label{Def-M-SVIEs} \rm A process $X(\cd)\in C_{\dbF}([0,T];L^2(\Om;\dbR^n))$ is called an
adapted solution to the forward equation in (\ref{FBSVIE-3-1}) if for every
$t\in[0,T]$, the corresponding equation is satisfied in the usual
It\^o's sense.
\ede

\mds

For the backward equation in (\ref{FBSVIE-3-1}), there are multi-type definitions of the solutions, see \cite{Lin 2002}, \cite{Shi-Wang 2012}, \cite{Shi-Wang-Yong 2015}, \cite{Wang-Yong 2015}, \cite{Wang-Zhang 2007}, etc. In this paper we would like to adopt the following:

\mds

\bde \label{Def-M-BSVIEs} \rm A pair of processes $(Y(\cd),Z(\cd\,,\cd))\in\cH^2(0,T;\dbR^m)$ is called a C-adapted solution to the second equation in (\ref{FBSVIE-3-1}) if for every $t\in[0,T]$,
the corresponding equation is satisfied for almost all $\omega\in\Omega$ and the  measurable process $\l(t,\cd)$ defined by
$$\ba{ll}
\ns\ds \l(t,\cd):= \psi(t,X(T))+\int_{\cd}^Tg(t,s,X(s),Y(s),Z(t,s),u(s))ds-\int_{\cd}^TZ(t,s)dW(s)
\ea$$
is in $L^2_{\dbF}(\Om;C(0,T;\dbR^m))$.
\ede

\mds

\begin{remark}
\em
Note that Definition \ref{Def-M-BSVIEs} shows some continuity of $(Y(\cd),Z(\cd,\cd))$ in some sense, hence we then name it the C-adapted solution. In contrast with the existed adapted solution under $\cH^2_{\D}(0,T;\dbR^m)$, more regularities of $(Y(\cd),Z(\cd,\cd))$ can be obtained by our new notion. For example, when the BSVIE in (\ref{FBSVIE-3-1}) degenerates into classical nonlinear BSDE, $Z(\cd,\cd)\in C([0,T];L^2_{\dbF}(0,T;\dbR^m))$ will become $Z(\cd)\in L^2_{\dbF}(0,T;\dbR^m)$ which obviously coincides with existed literature. However, this procedure does not work well for classical adapted solution where $Z(t,\cd)\in L^2_{\dbF}(\D;\dbR^m)$. In the second place, under our framework some terms, $Y(0)$, $Z(0,\cd)$ become meaningful and can be applied in the cost functional of optimal control problems.
\end{remark}

\mds

For FBSVIE (\ref{FBSVIE-3-1}), we introduce the following hypothesis.

\mds

{\bf(H1)} Let $\varphi(\cd)\in C_{\dbF}([0,T];L^2(\Omega;\dbR^n))$,
$\psi(\cd,0)\in C_{\dbF}([0,T];L^2(\Omega;\dbR^m))$,
$$\2n\ba{ll}
\ns\ds b,\si:[0,T]^2\times\dbR^n\times U\times\Omega\to\dbR^n,\ \
g:[0,T]^2\times\dbR^n\times\dbR^m\times\dbR^m\times
U\times\Omega\to\dbR^m,\ea $$
be measurable, for any $(t,x,u)\in [0,T]\times\dbR^n\times U$
$$s\mapsto\big(b(t,s,x,u),\si(t,s,x,u)),g(t,s,x,y,z,u)\big)$$
is $\cF$-progressively measurable on $[0,T]$,  for a.e. $(t,s,\omega)\in [0,T]^2\times\Omega$
$$(x,y,z,u)\mapsto(b(t,s,x,u),\si(t,s,x,u),g(t,s,x,y,z,u),\psi(t,x))$$
is continuously differentiable with uniformly bounded derivatives,
and, for
$$b_0(t,s):= b(t,s,0,0),\q\si_0(t,s):=\si(t,s,0,0),\q g_0(t,s):=
g(t,s,0,0,0,0),$$
one has
$$\sup_{t\in[0,T]}\dbE\[\(\int_0^T|b_0(t,s)|ds\)^2\!+\!\int_0^T|\si_0(t,s)|^2ds
\!+\!\(\int_0^T|g_0(t,s)|ds\)^2\]
<\infty.$$
Further, there exists a modulus of continuity
$\rho:[0,\infty)\to[0,\infty)$ such that
\begin{eqnarray*}
&&|b(t,s,x,u)-b(t',s,x,u)|+|\si(t,s,x,u)-\si(t',s,x,u)|+|\psi(t,x)-\psi(t',x)|\\
&&+|g(t,s,x,y,z,u)-g(t',s,x,y,z,u)|\le\rho(|t-t'|)[1+|x|+|y|+|z|+|u|],\\
&&\qq \forall\ t,t',s\in[0,T],\ \ \forall\ x\in\dbR^n,\ \forall\ y,z\in\dbR^m,\ \forall u\in U.
\end{eqnarray*}

\mds

In what follows, $C$ represents a generic positive constant, which may be different from line to line. Now let us discuss the well-posedness of (\ref{FBSVIE-3-1}). Some relevant study can also be found in \cite{Shi-Wang-Yong 2015}.

\mds

\bl\label{Pro-BSVIEs}
\rm {\rm(i)}  Let {\rm(H1)} hold. Then, for any $u(\cd)\in \cU_{ad}$, FBSVIE (\ref{FBSVIE-3-1})
admits a unique triple of
$(X(\cd),Y(\cd),Z(\cd,\cd))\in
C_{\dbF}([0,T];L^2(\Omega;\dbR^n))\times \cH^2(0,T;\dbR^m)$
in the spirt of Definition \ref{Def-M-SVIEs} and Definition \ref{Def-M-BSVIEs} such that for any $t\in[0,T],$
\bel{Lemma-3.1-estimate-X-Y-Z}~~~~\left\{\ba{ll}
\ns\ds \dbE|X(t)|^2 \!\le\! C\Big\{ \dbE|\f(t)|^2\! +\!
\dbE\Big(\int_0^t\!|b(t,s,0,u(s))|ds\Big)^2\! + \!\dbE
\int_0^t\!|\sigma(t,s,0,u(s))|^2ds\\
\ns\ds\qq\qq\q+\!\int_0^t\!\dbE|\f(s)|^2ds\! +\!
\int_0^t\dbE\Big(\int_0^s|b(s,r,0,u(r))|dr\Big)^2ds\\
\ns\ds\qq\qq\q+ \!\int_0^t\dbE
\int_0^s|\sigma(s,r,0,u(r))|^2drds \Big\},\\
\ns\ds  \dbE|Y(t)|^2\!+\! \dbE\int_t^T|Z(t,s)|^2ds\le C\Big\{
\dbE|\psi(t,X(T))|^2 \!+\!\int_t^T\dbE|\psi(s,X(T))|^2ds \\
\ns\ds\qq\qq\q+\!
\dbE\Big(\int_t^T\!|g(t,r,X(r),0,0,u(r))|dr\Big)^2\\
\ns\ds\qq\qq\q+\!
\int_t^T\!\dbE\Big(\int_s^T|g(s,r,X(r),0,0,u(r))|dr\Big)^2ds\Big\}.
\ea\right.\ee
\el
\rm

\begin{proof} \rm Note that the conclusions for forward SVIEs is obvious (see \cite{Shi-Wang-Yong 2015}), so next we will focus on the backward case.

\it Step 1: \rm  We prove the existence of $(Y_1(\cd),Z_1(\cd,\cd))\in \cH^2(0,T;\dbR^m)$ satisfying
\bel{Lemma-3.1-Y-1}\ba{ll}
\ns\ds Y_1(t)\!=\!\psi(t,X(T))\!+\!\!\int_t^T\!\!g(t,s,X(s),Y_0(s),Z_1(t,s),u(s))ds\!
-\!\!\int_t^T\!\!Z_1(t,s)dW(s),\ \hb{a.s.},
\ea\ee
where $(Y_0(\cd),Z_0(\cd,\cd))\in \cH^2_{\D}(0,T;\dbR^m)$ is the solution of BSVIE
\bel{Lemma-3.1-Y-0}\ba{ll}
\ns\ds Y_0(t)\!=\!\psi(t,X(T))\!+\!\!\int_t^T\!\!g(t,s,X(s),Y_0(s),Z_0(t,s),u(s))ds\!
-\!\!\int_t^T\!\!Z_0(t,s)dW(s).\ \hb{a.s.}\
\ea\ee
First, by Theorem 2.3 in  \cite{Shi-Wang-Yong 2015} (see also Proposition 3.2 in \cite{Wang-Yong 2015}), BSVIE (\ref{Lemma-3.1-Y-0}) admits a unique solution in $\cH^2_{\D}(0,T;\dbR^m)$.
Given $Y_0(\cd)$ and any $t\in[0,T]$, it is well known that the following parameterized BSDE,
$$\ba{ll}
\ns\ds \l_1(t,r)\!=\!\psi(t,X(T))\!+\!\!\int_r^T\!\!g(t,s,X(s),Y_0(s),Z_1(t,s),u(s))ds
\!-\!\!\int_r^T\!\!Z_1(t,s)dW(s)
\ea$$
with $r\in[0,T]$ admits a unique pair of measurable solution
$(\l_1(t,\cd),Z_1(t,\cd))\in L^2_{\dbF}(\Omega;C(0,T;\dbR^m))\times L^2_{\dbF}(0,T;\dbR^m)$
and the following estimate holds true,
\bel{Lemma-3.1-Y-1-estimate}\ba{ll}
\ns\ds \dbE\sup_{r\in[0,T]}|\l_1(t,r)|^2+\dbE\int_0^T|Z_1(t,s)|^2ds\\
\ns\ds \leq C\Big[\dbE|\psi(t,X(T))|^2+\dbE\Big(\int_0^T|g(t,s,X(s),Y_0(s),0,u(s))|ds\Big)^2\Big].
\ea\ee
In addition, by (H1) and the standard estimates for BSDEs, for any $t_0\in[0,T]$,
\bel{Lemma-3.1-lambda-1}\ba{ll}
\ns\ds \lim_{t\rightarrow t_0}\Big[\dbE\sup_{r\in[0,T]}|\l_1(t,r)-\l_1(t_0,r)|^2+
\dbE\int_0^T|Z_1(t,s)-Z_1(t_0,s)|^2ds\Big]=0.
\ea\ee
Next let $t=r$ and $Y_1(t):=\l_1(t,t)$ with any $t\in[0,T]$, we then have (\ref{Lemma-3.1-Y-1}) above. Since for any $t_0\in[0,T]$, one thus has $Y_1(\cd)\in C_{\dbF}([0,T];L^2(\Omega;\dbR^m))$ due to
$$\ba{ll}
\ns\ds \dbE|Y_1(t)-Y_1(t_0)|^2\leq \Big[\dbE\sup_{s\in[0,T]}|\l_1(t,s)-\l_1(t_0,s)|^2
+\dbE|\l_1(t_0,t)-\l_1(t_0,t_0)|^2\Big],
\ea$$
The result for $Z_1(\cd,\cd)$ comes from (\ref{Lemma-3.1-Y-1-estimate}) and (\ref{Lemma-3.1-lambda-1}).

\it Step 2: \rm In this step,  we prove the existence of the C-adapted solution $(Y_2(\cd),$ $Z_2(\cd,\cd))\in \cH^2(0,T;\dbR^m)$ in the sense of Definition \ref{Def-M-BSVIEs} of BSVIE
\bel{Lemma-3.1-Y-2}\ba{ll}
\ns\ds Y_2(t)\!=\!\psi(t,X(T))\!+\!\int_t^T\!\!g(t,s,X(s),Y_2(s),Z_2(t,s),u(s))ds\!-\!\int_t^T\!\!Z_2(t,s)dW(s),\ \ \hb{a.s.}
\ea\ee
with $t\in[0,T]$. First, let $(Y_1(\cd),Z_1(\cd,\cd))$ be the process defined in \it Step 1, \rm for any $t\in[0,T]$, the parameterized BSDE
$$\ba{ll}
\ns\ds \l_2(t,r)\!=\!\psi(t,X(T))\!+\!\!\int_r^T\!\!\!g(t,s,X(s),Y_1(s),Z_2(t,s),u(s))ds
\!-\!\!\int_r^T\!\!\!Z_2(t,s)dW(s),\ \hb{a.s.}
\ea$$
with  $ r\in[0,T]$ admits a unique pair of solution
$(\l_2(t,\cd),Z_2(t,\cd))\in L^2_{\dbF}(\Omega;C([0,T];\dbR^m))\times L^2_{\dbF}(0,T;\dbR^m).$
By (\ref{Lemma-3.1-Y-1}), (\ref{Lemma-3.1-Y-0}) and the estimates of adapted solution for BSVIE under $\cH^2_{\D}(0,T;\dbR^m)$, we obtain that
\bel{Lemma-3.1-Y-1-Y-0}\ba{ll}
\ns\ds \dbE\int_0^T|Y_1(t)-Y_0(t)|^2dt+\dbE\int_0^T\int_t^T|Z_1(t,s)-Z_0(t,s)|^2dsdt=0.
\ea\ee
We also need to compare $(\l_1,Z_1)$ with $(\l_2,Z_2)$. It follows from standard estimates of BSDEs and above (\ref{Lemma-3.1-Y-1-Y-0}) that
\bel{Lemma-3.1-lambda-1-lambda-2}\ba{ll}
\ns\ds \dbE\sup_{r\in[0,T]}|\l_1(t,r)-\l_2(t,r)|^2
+\dbE\int_0^T|Z_1(t,s)-Z_2(t,s)|^2ds\\
\ns\ds \leq C\dbE\Big[\int_0^T|g(t,s,X(s),Y_0(s),Z_1(t,s),u(s))\\
\ns\ds\qq
-g(t,s,X(s),Y_1(s),Z_1(t,s),u(s))|ds\Big]^2 =0, \ \forall\ t\in[0,T].
\ea\ee
Now, define $Y_2(t):=\l_2(t,t)$  for any $t\in[0,T]$. Notice that $Y_2(t)$ satisfies
$$\ba{ll}
\ns\ds Y_2(t)=\!\psi(t,X(T))\!+\!\!\int_t^T\!\!g(t,s,X(s),Y_1(s),Z_2(t,s),u(s))ds\!
-\!\!\int_t^T\!\!Z_2(t,s)dW(s).\ \hb{a.s.}\
\ea$$
This, together with (\ref{Lemma-3.1-lambda-1-lambda-2}) imply that
\bel{Lemma-3.1-Y-1-Y-2}\ba{ll}
\ns\ds \sup_{t\in[0,T]}\dbE|Y_1(t)-Y_2(t)|^2+\sup_{t\in[0,T]}\dbE\int_0^T|Z_1(t,s)-Z_2(t,s)|^2ds=0.
\ea\ee
As a result, $(Y_2,Z_2)\in \cH^2(0,T;\dbR^n)$. On the other hand, by  (\ref{Lemma-3.1-Y-1-Y-2}),
$$\ba{ll}
\ns\ds\q\dbE\Big|Y_2(t)\!-\!\psi(t,X(T))\!
-\!\!\int_t^T\!\!\!g(t,s,X(s),Y_2(s),Z_2(t,s),u(s))ds\!
+\!\!\int_t^T\!\!\!Z_2(t,s)dW(s)\Big|^2\\
\ns\ds=\dbE\Big| \int_t^T\!\!\big[g(t,s,X(s),Y_1(s),Z_2(t,s),u(s))\!
-\!g(t,s,X(s),Y_2(s),Z_2(t,s),u(s))\big]ds\Big|^2=0,
\ea$$
and
\begin{eqnarray*}
\l(t,r)&:=& \psi(t,X(T))\!+\!\int_r^T\!\!g(t,s,X(s),Y_2(s),Z_2(t,s),u(s))ds\!-\!\int_r^T\!\!Z_2(t,s)dW(s)\\
&=& \psi(t,X(T))\!+\!\int_r^T\!\!g(t,s,X(s),Y_1(s),Z_2(t,s),u(s))ds\!-\!\int_r^T\!\!Z_2(t,s)dW(s)\\
&=& \l_{2}(t,r),\ a.s.,\ \forall (t,r)\in [0,T]^2.
\end{eqnarray*}
Therefore, (\ref{Lemma-3.1-Y-2}) holds true a.s. for any $t\in[0,T]$ and  $\l(t,\cd)\in L^2_{\dbF}(\Omega;C([0,T];\dbR^m))$.

\it Step 3: \rm The uniqueness issue and related estimate.

Suppose there is another triple of $(Y_2',Z_2',\l_3')$. For suitable constant $\beta>0$, by using the tricks in Lemma 3.1 in \cite{Shi-Wang 2012} or Theorem 3.4 in \cite{Wang-Yong 2015} one has
$$\ba{ll}
\ns\ds \dbE\int_0^Te^{\beta s}|Y_2(s)-Y_2'(s)|^2ds+\dbE\int_0^Te^{\beta t}\int_t^T|Z_2(t,s)-Z_2'(t,s)|^2dsdt=0.
\ea$$
Considering $Y_2, Y_2'\in C_{\dbF}([0,T];L^2(\Omega;\dbR^m))$, one has $\dbP\big(\{\o,\ Y_2(t,\o)=Y_2'(t,\o)\}\big)=1$ for any $t\in[0,T].$ Given $(X(\cd),Y_2(\cd),u(\cd))$, obviously there exists
a unique pair of measurable processes
$(P(t,\cd),Q(t,\cd))\in L^2_{\dbF}(\Omega;C([0,T];\dbR^m))\times L^2_{\dbF}(0,T;\dbR^m)$
satisfying
\bel{Lemma-3.1-parametered-BSDE}\ba{ll}
\ns\ds~~~~~ P(t,r)\!=\!\psi(t,X(T))\!
+\!\!\int_r^T\!\!g(t,s,X(s),Y_2(s),Q(t,s),u(s))ds\!
-\!\!\int_r^T\!\!Q(t,s)dW(s),
\ea\ee
a.s., for any $r\in[0,T]$. Moreover, for any $s\in[0,T]$, the following estimate holds true,
\bel{Lemma-3.1-parametered-BSDE-estimate}\ba{ll}
\ns\ds \dbE\sup_{r\in[s,T]}|P(t,r)|^2+\dbE\int_s^T|Q_2(t,r)|^2dr\\
\ns\ds \leq C\Big[\dbE|\psi(t,X(T))|^2+\Big(\dbE\int_s^T|g(t,r,X(r),Y_2(r),0,u(r))|dr\Big)^2\Big].
\ea\ee
By the uniqueness of BSDE (\ref{Lemma-3.1-parametered-BSDE}) and $Y_2(\cd)=Y_2'(\cd)$, one finally has
$$\ba{ll}
\ns\ds
 \l_3(t,s)=\l_3'(t,s)=P(t,s),\ \forall s\in[0,T],\ \hb{a.s.}, \\
\ns\ds Z_2(t,s)=Z_2'(t,s)=Q(t,s),\  a.e.\ s\in[0,T],\   \ \hb{a.s.}, \ \forall t\in[0,T].
\ea
$$
On the other hand, let $s=t$, it follows from (\ref{Lemma-3.1-parametered-BSDE-estimate}) that
$$\ba{ll}
\ns\ds \dbE |Y_2(t)|^2+\dbE\int_t^T|Z_2(t,r)|^2dr\\
\ns\ds \leq C\Big[\dbE|\psi(t,X(T))|^2+\Big(\dbE\int_t^T|g(t,r,X(r),0,0,u(r))|dr\Big)^2+\dbE\int_t^T|Y_2(r)|^2dr\Big].
\ea$$
Then the second estimate of (\ref{Lemma-3.1-estimate-X-Y-Z}) associated with $(Y_2,Z_2)$ follows from the Gronwall inequality. This completes the proof of Lemma \ref{Pro-BSVIEs}.
\end{proof}

\mds

Given $(X,Y,Z)$ satisfies FBSVIE (\ref{FBSVIE-3-1}) in the spirt of Definition \ref{Def-M-SVIEs} and Definition \ref{Def-M-BSVIEs}, let us introduce the cost functional as
follows:
\bel{cost3.1}J(u(\cd))=\dbE \int_0^T
f(s,X(s),Y(s),Z(0,s),u(s))ds + \dbE h(X(T),Y(0)).\ee
For the involved functions $f$, $h$ in (\ref{cost3.1}), we impose the
following hypothesis.

\mds

{\bf (H2)} Let
$ f:[0,T]\times\dbR^n\times\dbR^m\times\dbR^m \times U
\times\Omega\to\dbR$,
$
h:\dbR^n\times\dbR^m\times\Omega\to\dbR$
be measurable such that
$ (x,y,z,u)\mapsto f(s,x,y,z,u)$,
$(x,y)\mapsto h(x,y)$
are continuously differentiable with the derivatives being bounded
by $L(1+|x|+|y|+|z|+|u|)$, $L>0$.

\mds

Now, we state our optimal control problem.

\mds

\bf Problem (C). \rm With the state equation (\ref{FBSVIE-3-1}),
find $\bar u(\cd)$ such that
\bel{3.4}J(\bar u(\cd))=\inf_{u(\cd)\in\cU_{ad}}J(u(\cd)).\ee

For any given optimal 4-tuple $(\bar X(\cd),\bar Y(\cd),\bar
Z(\cd\,,\cd),\bar u(\cd))$ of Problem (C), $t,\ s\in[0,T]$, we denote
$$\ba{ll}
\ns\ds b_x(t,s):= b_x(t,s,\bar X(s),\bar u(s)),\q b_u(t,s):= b(t,s,\bar
X(s),\bar u(s)),\\
\ns\ds \psi_x(t):=\psi_x(t,\bar X(T)), \qq\qq\  h_x:= h_x(\bar X(T),\bar Y(0)),\\
\ns\ds f_x(s):= f_x(s,\bar X(s),\bar Y(s),\bar Z(0,s),\bar u(s)).
\ea$$
The notations $\si_x(t,s)$, $\si_u(t,s)$, etc are similarly defined.
Also, for any scalar valued function, say, $f_x(s)$ is
regarded as a column vector.
Let
$(\bar X(\cd),\bar Y(\cd),\bar Z(\cd\,,\cd),\bar u(\cd))$ be an optimal
4-tuple of Problem {\rm(C)}, $(\l(\cdot), \xi(\cdot), \mu(\cdot), $ $\nu(\cdot), p(\cdot), q(\cdot,\cdot))$ be the solution to the following first order adjoint system with respect to the control FBSVIE (\ref{FBSVIE-3-1}).
\bel{adjoint-1}\left\{\2n\ba{ll}
\ns\ds\l(0)\!=\!\dbE(h_y)+\int_0^T\! f_z(s)dW(s)\!+\!\int_0^T\!g_z(0,s)^{\top}\dbE_s\l(0)dW(s),\\
\ns\ds\xi(t)\!=\!g_y(0,t)^{\top}\dbE_t\l(0)\!+\!f_y(t)\!+\!\int_0^t\!g_y(s,t)^{\top}\dbE_t[\xi(s)]ds
\!+\!\int_t^T\!g_z(t,s)^{\top}\dbE_s[\xi(t)]dW(s),\\
\ns\ds\m(t)\!=\!h_x\!+\!\psi_x(0)^{\top}\l(0)\!+\!\int_0^T\!\psi_x(s)^{\top}\xi(s)ds
\!-\!\int_t^T\n(s)dW(s),
\\
\ns\ds p(t)\!=\!b_x(T,t)^{\top}\m(T)\!+\!\si_x(T,t)^{\top}\n(t)\!+\!g_x(0,t)^{\top}\dbE_t\l(0)
\!+\!f_x(t)
\!+\!\int_0^tg_x(s,t)^{\top}\xi(s)ds\\
\ns\ds\qq
\!+\!\int_t^T\(b_x(s,t)^{\top}p(s)
\!+\!\si_x(s,t)^{\top}q(s,t)\)ds
\!-\!\int_t^T\!q(t,s)dW(s),\ea\right.\ee
where $\dbE_s\eta$ stands for $\dbE[\eta|\cF_{s}]$, and $A^{\top}$ is the transpose of vector/matrix $A$. Note that above $\dbE_{\cd}\l(0)$ solves an FSDE,
$\xi(\cd)$ is allowed to be $\cF_T$-measurable. $(\m(\cd),\n(\cd))$, $(p(\cd),q(\cd\,,\cd))$
satisfies a BSDE and BSVIE respectively. In the following, let us define,
\bel{cH-1}\qq\ba{ll}
\ns\ds\mH_{u}(s;\bar X(\cd),\bar Y(\cd),\bar Z(\cd,\cd),\bar
u(\cd))\!:= \!\!g_u(0,s)^{\top}\dbE_s\l(0)\!
+\!\dbE_s\int_0^sg_u(t,s)^{\top}\xi(t)dt\!+\!f_u(s)\\
\ns\ds \!+\!\dbE_s\int_s^T\!\!\big[b_u(t,s)^{\top}p(t)\!+\!\si_u(t,s)^{\top}
q(t,s)\big]dt\! +\!b_u(T,s)^{\top}\dbE_s\m(T)\!+\!\si_u(T,s)^{\top}\n(s).
\ea\ee

\mds

\begin{theorem}\label{Theorem-main-result}
 \rm Let {\rm(H1)}, {\rm(H2)} hold and $(\bar
X(\cd),\bar Y(\cd),\bar Z(\cd\,,\cd),\bar u(\cd))$ be an optimal
4-tuple of Problem {\rm(C)}. Then for almost every $s\in[0,T]$,
\bel{maximum-1}\ba{ll}
\ns\ds\lan  \mH_{u}(s;\bar X(\cd),\bar Y(\cd),\bar Z(\cd,\cd),\bar
u(\cd)) ,v
\ran\ge0, \qq \forall  v\in \T_{U}(\bu(s)),\ \ \hb{a.s.}
\ea\ee
\end{theorem}
\rm

The proof of Theorem \ref{Theorem-main-result} will be given in Section 4.

\begin{example}
\rm When the FBSVIE (\ref{FBSVIE-3-1}) degenerates
the classical FBSDE  (\ref{State-FBSDEs}), the adjoint system (\ref{adjoint-1}) becomes
\bel{adjoint-1-FBSDEs}\q~~~\left\{\2n\ba{ll}
\ns\ds\l(0)\!=\!\dbE(h_y^{\top})\!+\!\int_0^T\!\! f_z(s)dW(s)\!+\!\int_0^T\!\!g_z(s)^{\top}\dbE_s\l(0)dW(s),\\
\ns\ds\xi(t)\!=\!g_y(t)^{\top}\dbE_t\l(0)\!+\!f_y(t)^{\top}\!
+\!g_y(t)^{\top}\!\!\int_0^t\!\!\dbE_t[\xi(s)]ds
\!+\!\int_t^T\!\!\!g_z(s)^{\top}\dbE_s[\xi(t)]dW(s),\\
\ns\ds\m(t)\!=\!h_x\!+\!\psi_x ^{\top}\l(0)\!+\!\psi_x ^{\top}\int_0^T\!\!\xi(s)ds
\!-\!\int_t^T\!\!\n(s)dW(s),
\\
\ns\ds p(t)\!=\!b_x(t)^{\top}\Big[\m(T)\!+\!\int_t^T\!\!p(s)ds\Big]
\!+\!\si_x(t)^{\top}\Big[\n(t)\!+\!\int_t^T\!\!q(s,t)ds\Big]\!+\!f_x(t)\\
\ns\ds\qq \ \ +\!g_x(t)^{\top}\dbE_t\l(0)\!+\!g_x(t)^{\top}\int_0^t\!\!\xi(s)ds
\!-\!\int_t^T\!q(t,s)dW(s).\ea\right.\ee

and the first order necessary condition (\ref{maximum-1}) reduces to
\bel{maximum-1-FBSDEs}\q\ba{ll}
\ns\ds\Big\langle b_u(s)^{\top}\Big[\dbE_s\int_s^{T}\!\!p(t)dt\!+\!\dbE_s\m(T)\Big]
\!+\!\si_u(s)^{\top}
\Big[\int_t^T\!\!q(t,s)dt\!+\!\n(s)\Big]\!+\!f_u(s)\\
\ns\ds\!+\!g_u(s)^{\top}\Big[\dbE_s\l(0)\!+\!\dbE_s\int_0^s\!\! \xi(t)dt\Big],v
\Big\rangle\!\ge\!0, \ \forall  v\!\in\!\T_{U}(\bu(s)),\ \hb{a.e.}\ s\!\in\! [0,T],\ \hb{a.s.}
\ea\ee
For $\l(0)$ and $\xi(\cd)$ in the first two equations of (\ref{adjoint-1-FBSDEs}), if we denote by
$$
\L(t):= \dbE_t\l(0)+\dbE_t\int_0^t\xi(s)ds,\qq t\in[0,T],
$$
it is then a direct calculation that for any $t\in[0,T]$,
\bel{Example-FBSDE-adjoint-equation-1}\ba{ll}
\ns\ds \L(t)\!=\!\dbE(h_y^{\top})\!+\!\int_0^t\!\!\big[f_y(s)\!+\!g_y(s)^{\top}\L(s)\big]ds\!+\!
\int_0^t\!\!\big[f_z(s)\!+\!g_z(s)^{\top}\L(s)\big]dW(s).
\ea\ee
And, the third equation of (\ref{adjoint-1-FBSDEs}) becomes
$$\ba{ll}
\ns\ds \mu(t)=h_x+\psi_x^{\top}\L(T)-\int_t^T\nu(s)dW(s),\qq t\in[0,T].
\ea$$
For any $t\in[0,T],$ if we define
$$\ba{ll}
\ns\ds P(t):=\dbE\Big[\m(T)+\int_t^Tp(s)ds\Big],\ \  Q(t):=\Big[\n(t)+\int_t^Tq(s,t)ds\Big],
\ea$$
hence one can see that $(P(\cd),Q(\cd))$ satisfies a BSDE of
$$\ba{ll}
\ns\ds P(t)= \m(T)+\int_t^Tp(s)ds-\int_t^TQ(s)dW(s),\ \ t\in[0,T].
\ea$$
On the other hand, recalling previous defined $(P(\cd),Q(\cd))$, let us take condition expectation on both sides of the fourth equation in (\ref{adjoint-1-FBSDEs}),
$$\ba{ll}
\ns\ds p(t)=b_x(t)^{\top}P(t)+\sigma_x(t)^{\top}Q(t)+f_x(t)+g_x(t)^{\top}\L(t).
\ea$$
As a result,
\bel{Example-FBSDE-adjoint-equation-2}\ba{ll}
\ns\ds P(t)=h_x+\psi_x^{\top}\L(T)+\int_t^T\big[g_x(s)^{\top}\L(s)+f_x(s)^{\top}\big]ds\\
\ns\ds\qq\qq +\int_t^T\big[b_x(s)^{\top}P(s)+\sigma_x(s)^{\top}Q(s)\big]ds-\int_t^TQ(s)dW(s),
\ea\ee
and the maximum condition (\ref{maximum-1-FBSDEs}) can be rewritten as,
\bel{Example-FBSDE-necessary-condition}\q~~\ba{ll}
\ns\ds \lan b_u(s)^{\top}P(s)\!+\!\sigma_u(s)^{\top}Q(s)\!+\!f_u(s)\!+\!g_u(s)^{\top}\L(s),v\ran \!\geq\!0,\ \forall\  v\!\in \!\T_{U}(\bu(s)),\  \hb{a.s.}
\ea\ee
Compared with the references \cite{Wu 2013}, \cite{Yong 2010} with non-convex control region $U$, here we can obtain a slight different condition via only first-order adjoint equation (\ref{Example-FBSDE-adjoint-equation-1}) and (\ref{Example-FBSDE-adjoint-equation-2}) and less differentiability requirements on the coefficients. As a trade-off, our condition (\ref{Example-FBSDE-necessary-condition}) is weaker than the existed results in the sense that it can be implied by the later.
\end{example}

\begin{remark}\label{strucural remark}
\em
When $\T_{U}(\bu(s))=\{0\}$ for a.e. $(s,\omega)\in [0,T]\times\Omega$, condition (\ref{maximum-1}) turns out to be trivial. For example, if $U=\{0,1\}$, by the definition of the adjacent cone one has $\T_{U}(0)=\T_{U}(1)=\{0\}$. In this case, we have to admit that our variation approach is infeasible. However, by imposing some structural assumptions on $U$, we can ensure the effectiveness of (\ref{maximum-1}), i.e. the existence of set $A\subset[0,T]\times\Omega $ with positive measure such that
\begin{equation}\label{strucural assum}
\T_{U}(\bu(s))\neq\{0\},\q (s,\omega)\in A.
\end{equation}
We refer the reader the Chapter 4 in \cite{Aubin90} for a detailed discussion in this respect. In what follows, we introduce an important example which is wildly used in practice.
\end{remark}

\begin{example}\label{Example-adajacent-cone-1}
\em
Let $g^{i}:\mr^{l}\to \mr$, $i=1,2,\cdots, k$ be continuously differentiable functions. Define
$$U:=\big\{u\in \mr^l\ \big|\ g^{i}(u)\le 0, \ i=1,2,\cdots, k\big\}.$$
Obviously, $U$ is a closed set of $\mr^l$, and is not necessary to be convex.
Assume that the linearly independent constraint qualification (LICQ, in short) is satisfied at any $u\in bdU$ (the boundary of $U$), that is,
$$
g^{i}_{u}(u), \quad i\in I[u]=\big\{i\ \big|\ g^{i}(u)=0,\ i=1,2,\cdot\cdot\cdot, k\big\}
$$
are linearly independent, then for any fixed $u\in bdU$ there exists a vector $v^{u}\in \mr^l$ such that
$$\inner{g_{u}^{i}(u)}{v^{u}}<0,\quad \forall\ i\in I[u].$$
Combining Theorem 4.3.3 and Proposition 4.3.7 in \cite{Aubin90}, we obtain that
$$\T_{U}(\bar{u}(s))=\big\{v\in \mr^l\ |\ \inner{g^{i}_{u}(\bar{u}(s))}{v}\le 0, \ \forall\ i\in  I[\bar{u}(s)]\big\},\quad \ \hb{a.e.}\ s\in[0,T],\ \hb{a.s.}$$
In this case, condition (\ref{maximum-1}) implies that
\bel{inclusion condition}\ba{ll}
\ns\ds 0\in \mH_{u}(s;\bar X(\cd),\bar Y(\cd),\bar Z(\cd,\cd),\bar
u(\cd))+N_{U}(\bar{u}(s))
,\ \hb{a.e.}\ s\in[0,T], \ \hb{a.s.}
\ea\ee
where
$$N_{U}(\bar{u}(s))=\Big\{\sum_{i=1}^{k}\lambda_{i}g_{u}^{i}(\bar{u}(s))\ \Big|\
\lambda_i\ge 0,\ i=1,2,\cdot\cdot\cdot,k,\ \lambda_i=0,\ \hb{with} \ i\in I[\bar{u}(s)]\Big\}.$$
By Filippov's selection lemma, we conclude that there exists an $\mmf$-adapted process $\lambda(\cdot)=(\lambda_1(\cdot),\lambda_2(\cdot),~$ $\cdot\cdot\cdot,\lambda_k(\cdot))$ satisfying the following type condition:
\bel{maximum-2add}\qq\left\{\ba{ll}
\ns\ds
\lambda_i(s)\ge 0,\ \ \lambda_i(s)g^i(\bar{u}(s))=0,\ \ i=1,2,\cdots, k,\ \ \hb{a.e.}\ s\in[0,T],\ \ \hb{a.s.},\\
\ns\ds \mH_{u}(s;\bar X(\cd),\bar Y(\cd),\bar Z(\cd,\cd),\bar
u(\cd))\!+\!\sum_{i=1}^{k}\lambda_{i}(s)g_{u}^{i}(\bar{u}(s))
\!=\!0,\ \ \hb{a.e.}\ s\in[0,T],\ \hb{a.s.}
\ea\right.\ee
\end{example}

\mds

To get more feelings of LICQ and adjacent cone in Example \ref{Example-adajacent-cone-1}, let us look at,

\mds

\begin{example}
\em
Suppose $U\subset\mathbb{R}^2$ is a torus defined by:
$$\left\{\ba{ll}
\ns\ds g^{1}(u_{1},u_{2})=u_{1}^{2}+u_{2}^{2}-4\le 0,\\
\ns\ds g^{2}(u_{1},u_{2})=2-u_{1}^{2}-u_{2}^{2}\le 0.
\ea\right.$$
Obviously, $U$ is nonconvex and for any $u=( u_{1}, u_{2})\in bd U$, one has either $ u_{1}^{2}+ u_{2}^{2}-4=0$ or $2- u_{1}^{2}-  u_{2}^{2}=0$.
If $ u_{1}^{2}+ u_{2}^{2}-4=0$, then $I[ u]=\{1\}$ and $g^{1}_{ u}=(2 u_{1},2 u_{2})\neq 0$; while if $2-  u_{1}^{2}- u_{2}^{2}=0$, one then has $I[u]=\{2\}$ and $g^{1}_{  u}=(-2  u_{1},-2  u_{2})\neq 0$. Consequently, in both cases, the LICQ is satisfied and
$$\T_{U}(  u)=\left\{\ba{ll}
\ns\ds \mathbb{R}^2,\qq\qq\qq\qq\qq\qq\qq     1< u_{1}^{2}+  u_{2}^{2}<4,\\
\ns\ds \{(v_{1},v_{2})\in \mathbb{R}^2\ |\   u_{1}v_{1}+  u_{2}v_{2}\le 0\}, \q
  u_{1}^{2}+  u_{2}^{2}-4=0,\\
\ns\ds \{(v_{1},v_{2})\in \mathbb{R}^2\ |\   u_{1}v_{1}+ u_{2}v_{2}\ge 0\},  \q
2- u_{1}^{2}-  u_{2}^{2}=0.
\ea\right.$$
\end{example}

\section{Proofs of the main results}

This part is devoted to proving Theorem \ref{Theorem-main-result}. We
make some preparations in the first two subsections and present the proof of Theorem \ref{Theorem-main-result} in the last part.

\subsection{A general duality principle}

In this part, we would like to establish a unified duality principle
for stochastic Volterra integral systems, which covers the forward
case in \cite{Yong 2008}, and the forward-backward case in
\cite{Shi-Wang-Yong 2015}. To do so, let us firstly look at the
following kind of equation
\bel{SFVIE-1}~\ba{ll}
\ns\ds\xi(t)=\a(t)+\int_0^t\!A(s,t)^{\top}\dbE_t[\xi(s)]ds
+\int_0^t\dbE_s[B(s,t)^{\top}\xi(s)]dW(s)\\
\ns\ds\qq\q+\int_0^T\!\b(t,s)dW(s)+\int_t^T\!D(t,s)^{\top}\dbE_s[\xi(t)]dW(s),
\ \hb{a.e.}\ t\in[0,T],\ \hb{a.s.},\ea\ee
where $(\a(\cd),\b(\cd\,,\cd))\in L^2_\dbF(0,T;\dbR^m)\times
L^2(0,T;L^2_{\dbF}(0,T;\dbR^m))$. Such equation
was introduced in \cite{Shi-Wang-Yong 2015} and its wellposedness in
space $L^2_{\cF_T}(0,T;\dbR^m)$ was discussed. However, in order to
obtain a general duality principle aforementioned, we need more higher
regularities of the $non$-$adapted$ solution $\xi(\cd)$. More
precisely, next we want to seek for solution $\xi(\cd)\in C
_{\cF_T}([0,T];L^2(\Omega;\dbR^m))$ under appropriate conditions.

\mds

{\bf(H3)} $A,\ B,\ D:[0,T]^2\times\Omega\to\dbR^{m\times m}$
are measurable and bounded processes such that for almost all
$t\in[0,T]$, $s\1n\mapsto\1n(A(t,s),\1n B(t,s),\1n D(t,s))$ is
$\dbF$-adapted on $[0,T]$. Moreover, for almost all $s\in[0,T]$,
$A(s,\cd),\ B(s,\cd),\ D(\cd,s)$ are continuous processes.

\mds

\bt \label{SFVIE-regularity}
 \rm Let {\rm(H3)} hold. Then for any
$(\a(\cd),\b(\cd\,,\cd))\!\in \!C _\dbF([0,T];L^2(\Omega;\dbR^m))$ $\times C
([0,T];L^2_{\dbF}(0,T;\dbR^m))$, (\ref{SFVIE-1}) admits a unique
solution $\xi(\cd)\in C _{\cF_T}([0,T];L^2(\Omega;\dbR^m))$.
\et

\begin{proof} \rm To make the argument more readable, we will separate
the procedures into several parts.

\it Step 1: \rm By Theorem 4.1 in \cite{Shi-Wang-Yong 2015},
equation (\ref{SFVIE-1}) admits a unique solution $\xi(\cd)\in
L^2_{\cF_T}(0,T;\dbR^m)$. Given such $\xi(\cd)$, for any $t\in[0,T]$, let
\bel{4.2}\ba{ll}
\ns\ds f(t;\xi):= \alpha(t)\!+\!\int_0^t\!\!A
(s,t)^{\top}\dbE_t\xi(s)ds\!+\!\int_0^t\!\!\Big[\beta(t,s)
\!+\!\dbE_s\big[B(s,t)^{\top}\xi(s)\big]\Big]dW(s).
\ea\ee
In this step we will show that $f(\cd;\xi)\in C
_{\dbF}([0,T];L^2(\Omega;\dbR^m))$. Obviously, $f(t;\xi)$ is $\cF_t$-measurable with $t\in[0,T]$. By the
requirements of $\alpha(\cd)$, $\beta(\cd,\cd)$, as well as
{\rm(H3)}, one has $\sup\limits_{t\in[0,T]}\dbE|f(t;\xi)|^2<\infty.
$
We only need to
prove that
\bel{SFVIE-regularity-first-one}\ba{ll}
\ns\ds \dbE\big|f(t;\xi)-f(t_0;\xi)\big|^2\rightarrow0,\ \ t\rightarrow t_0,\ \
\forall t_{0}\in [0,T].
\ea\ee
To this end, let us firstly look at the case
with $t_0\in[0,T),$ $t\geq t_0$. By the definition of
$\beta(\cd,\cd)$ and dominated convergence theorem, one has
\bel{beta-1}~~\ba{ll}
\ns\ds \dbE\Big|\int_0^{t}
\beta(t,s)dW(s)-\int_{0}^{t_0}\beta(t_0,s)dW(s)\Big|^2
\\
\ns\ds\leq C\Big[\dbE
\int_0^{t}\big[\beta(t,s)-\beta(t_0,s)\big]^2ds +
\dbE\int_{0}^T|\beta(t_0,s)|^2I_{[t_0,t]}ds\Big]\rightarrow0,\ \
t\to t_0^{+}.
\ea\ee
Similar as above, we also have,
\bel{C-s-t-1}\ba{ll}
\ns\ds \dbE\Big|\int_0^t \dbE_s\big[B(s,t)^{\top}\xi(s)\big]dW(s)
-\int_0^{t_0} \dbE_s\big[B(s,t_0)^{\top}\xi(s)\big]dW(s)\Big|^2\\
\ns\ds\leq C\dbE\int_0^t\!\big|\dbE_s
[B(s,t)^{\top}\xi(s)\!-\!B(s,t_0)^{\top}\xi(s)]\big|^2ds
+C\dbE\int_0^T\!\!\!\big|\dbE_s[B(s,t_0)^{\top}\xi(s)]\big|^2I_{[t_0,t]}ds\\
\ns\ds\q\rightarrow0,
\ \ t\to t_0^{+}.
\ea\ee
Since $\xi(\cd)\in L^2_{\cF_T}(0,T;\dbR^m)$, there exists a unique $\phi_{\xi}(\cd,\cd)\in
L^2(0,T;L^2_{\dbF}(0,T;\dbR^m))$ s.t.
$$\ba{ll}
\ns\ds \xi(s)=\dbE\xi(s)+\int_0^T\phi_{\xi}(s,r)dW(r),\ \ae s\in[0,T],\ \hb{a.s.}
\ea$$
Hence
$$
\dbE_t\xi(s)\!-\!\dbE_{t_0}\xi(s)\!=\!\int_{t_0}^t\phi_{\xi}(s,r)dW(r),\ \ \dbE\int_0^T\!\!\int_0^T\!|\phi_{\xi}(s,r)|^2drds\!\leq\! C\dbE\int_0^T\!|\xi(s)|^2ds.$$

As a result, the following estimates hold:
\bel{A-s-t-2}\left\{\ba{ll}
\ns\ds
\dbE\Big[\int_0^t|A(s,t)^{\top}|\big|\dbE_t\xi(s)-\dbE_{t_0}\xi(s)
\big|ds\Big]^2 \leq
C\dbE\int_{t_0}^t\int_0^t|\phi_{\xi}(s,r)|^2dsdr,\\
\ns\ds
\dbE\Big[\int_0^t|A(s,t)^{\top}-A(s,t_0)^{\top}||\dbE_{t_0}\xi(s)|ds\Big]^2\\
\ns\ds\qq\leq
C\dbE \int_0^t|A(s,t)^{\top}-A(s,t_0)^{\top}|^2|\dbE_{t_0}\xi(s)|^2ds,\\
\ns\ds
\dbE\int_0^T|A(s,t_0)^{\top}|^2|\dbE_{t_0}\xi(s)|^2I_{[t_0,t]}(s)ds\leq
C\dbE\int_0^T |\xi(s)|^2I_{[t_0,t]}(s)ds.
\ea\right.\ee
As a result, by the estimates in (\ref{A-s-t-2}),
\bel{A-s-t-1}\ba{ll}
\ns\ds \dbE\Big|\int_0^tA(s,t)^{\top}\dbE_t\xi(s)ds-\int_0^{t_0}
A(s,t_0)^{\top}\dbE_{t_0}\xi(s)ds\Big|^2\\[+0.3em]
\ns\ds\leq C
\dbE\Big[\!\int_0^t\!|A(s,t)^{\top}|\big|\dbE_t\xi(s)\!-\!\dbE_{t_0}\xi(s)\big|ds\Big]^2\!\!+\!
C\dbE\Big[\!\int_0^t\!\!|A(s,t)^{\top}\!-\!A(s,t_0)^{\top}||\dbE_{t_0}\xi(s)|ds\Big]^2\\[+0.3em]
\ns\ds\q+C\dbE\int_0^T|A(s,t_0)^{\top}|^2|\dbE_{t_0}\xi(s)|^2I_{[t_0,t]}(s)ds\rightarrow 0,\ \ \ t
\to
t_0^{+},
\ea\ee
Using (\ref{beta-1}), (\ref{C-s-t-1}), (\ref{A-s-t-1}) and
$\alpha(\cd)\in C _{\dbF}([0,T];L^2(\Omega;\dbR^m))$, one has
$$\dbE\big|f(t;\xi)-f(t_0;\xi)\big|^2\to0,\ \ t\to t_0^+.$$
Similarly, one can deal with the case as
$t\to t_0^-$. Therefore, one has above (\ref{SFVIE-regularity-first-one}) and thus $f(\cd,\xi)\in C _{\dbF}([0,T];L^2(\Omega;\dbR^m))$.

\it Step 2: \rm Let $f(\cd;\xi)$ be the process defined by (\ref{4.2}). For any $t\in[0,T]$, let us
consider the SDE
\bel{SDE-t-1}\ba{ll}
\ns\ds
\L(s,t)=f(t;\xi)+\int_t^s\big[D(t,r)^{\top}\L(r,t)+\beta(t,r)\big]dW(r),\qq
\forall s\in[t,T],
\ea\ee
which obviously admits a unique solution $\L(\cd,t)\in
L^2_{\dbF}(\Omega;C([t,T];\dbR^m))$. Let $s=T$ and $\eta(t):=\L(T,t)$. Then, for any $r\geq t$, $\L(r,t)=\dbE_r\eta(t)$,  and
\bel{SDE-t-5}\ba{ll}
\ns\ds
\eta(t)=f(t;\xi)+\int_t^T\big[D(t,r)^T\dbE_r\eta(t)+\beta(t,r)\big]dW(r),\qq
\forall t\in[0,T].
\ea\ee
In this step, we would like to show that $\eta(\cd)\in
C_{\cF_T}([0,T];L^2(\Omega;\dbR^m))$.

As to (\ref{SDE-t-1}), by BDG inequality, for any $\t\in[t,T]$, one
has,
$$
\dbE\sup_{s\in[t,\t]}|\L(s,t)|^2 \leq
C\dbE|f(t;\xi)|^2+C\dbE\int_t^T|\beta(t,r)|^2dr
+C\dbE\int_t^{\t}|D(t,r)^{\top}\L(r,t)|^2dr.
$$
It then
follows from Gronwall inequality that,
$$ \dbE\sup_{s\in[t,\t]}|\L(s,t)|^2 \leq
C\dbE|f(t;\xi)|^2+C\dbE\int_t^T|\beta(t,r)|^2dr<\infty.
$$
Let $\t=T$, one has,
\bel{SDE-t-2}\ba{ll}
\ns\ds   \dbE\sup_{s\in[t,T]}|\L(s,t)|^2\leq C \dbE|f(t;\xi)|^2+C
\dbE\int_t^T|\beta(t,s)|^2ds<\infty,\ \ \forall t\in[0,T].
\ea\ee
Next, let $t_0\in[0,T]$, $t\leq t_0$, and denote $\wt\L(\cd;t,t_0):=
\L(\cd,t)-\L(\cd,t_0)$. One has,
\begin{eqnarray*}
\wt\L(s;t,t_0)&=&
f(t;\xi)-f(t_0;\xi)+\int_{t}^{t_0}\big[D(t,r)^{\top}\L(r,t_0)+\beta(t,r)\big]dW(r)\\
&&+\int_t^s\big[
D(t,r)^{\top}\wt\L(r,t,t_0)+\f(r,t,t_0)I_{[t_0,T]}(r)\big]dW(r),\ \
s\in[t,T],
\end{eqnarray*}
where for $r\in[t_0,T]$,
$$ \f(r,t,t_0):= \big[D(t,r)^{\top}-D(t_0,r)^{\top}\big]\L(r,t_0)
+\big[\beta(t,r)-\beta(t_0,r)\big].
$$
Similar to (\ref{SDE-t-2}), we obtain that
\begin{eqnarray*}
&&\dbE\sup_{s\in[t,T]}|\wt \L(s,t,t_0)|^2\leq
C\dbE|f(t;\xi)-f(t_0;\xi)|^2 +\dbE\int_{t_0}^{T}|\f(r,t,t_0)|^2dr\\
&&\qq\qq\qq\qq\qq+C\dbE\int_{t}^{t_0}|D(t,r)^{\top}\L(r,t_0)+\beta(t,r)|^2dr.
\end{eqnarray*}
Then, by the dominated convergence theorem, one concludes that
$$
\dbE|\eta(t)-\eta(t_0)|^2=\dbE|\L(T,t)-\L(T,t_0)|^2\rightarrow 0,\qq
t\to t_0^{-}.
$$

Similar conclusion also holds true with $t\geq t_0$. Considering
(\ref{SDE-t-2}), we have $\eta(\cd)\in
C_{\cF_T}([0,T];L^2(\Omega;\dbR^m))$.

\it Step 3: \rm In this part, we would like to prove the existence
and uniqueness of the solution of equation (\ref{SFVIE-1}).

Let $\xi(\cd)$ be a solution of equation (\ref{SFVIE-1}) in $L^2_{\cF_T}(0,T;\dbR^m)$. By the definition of $f(t;\xi)$,
$\xi(\cd)$ is also a solution of equation (\ref{SDE-t-5}). On the other hand, (\ref{SDE-t-5}) is a
particular case of (\ref{SFVIE-1}) with $A(\cd,\cd)=B(\cd,\cd)=0$ and $\beta(t,s)=0$ with $t\geq s$, hence the uniqueness of adapted solution in $L^2_{\cF_T}(0,T;\dbR^m)$ of (\ref{SDE-t-5}) can be implied by the case of (\ref{SFVIE-1}). Consequently,
$$\ba{ll}
\ns\ds \dbE\int_0^T|\xi(s)-\eta(s)|^2ds=0,\qq
\dbE|f(t;\xi)-f(t;\eta)|^2=0,\ \ \forall t\in[0,T].
\ea$$
So $\eta(\cd)\in C_{\cF_T}([0,T];L^2(\Omega;\dbR^m))$ also satisfies
$$\ba{ll}
\ns\ds
\eta(t)=f(t;\eta)+\int_t^T\big[D^{\top}(t,r)\dbE_r\eta(t)+\beta(t,r)\big]dW(r),\
\ \hb{a.s.} \ \ \forall t\in[0,T],
\ea$$
which means $\eta(\cd)\in C_{\cF_T}([0,T];L^2(\Omega;\dbR^m))$ is a
solution of equation (\ref{SFVIE-1}). As to the uniqueness, suppose there are two solutions $\eta_i(\cd)\in C_{\cF_T}([0,T];L^2(\Omega;\dbR^m))$, $i=1,2$ satisfying (\ref{SFVIE-1}). By the uniqueness under $L^2_{\cF_T}(0,T;\dbR^m)$, as well as the continuity of $\eta_i(\cd)$ from $[0,T]$ to $L^2_{\cF_T}(\Omega;\dbR^m)$, one can obtain,
$$\ba{ll}
\ns\ds
\dbE\int_0^T|\eta_1(s)-\eta_2(s)|^2ds=0 \Rightarrow \dbE|\eta_1(t)-\eta_2(t)|^2=0,\ \ \forall t\in[0,T].
\ea$$
This completes the proof of Theorem \ref{SFVIE-regularity}. \end{proof}

\mds

Now let us look at the following backward equations,
\bel{BSVIE-adjoint-1}\left\{\ba{ll}
\ns\ds
Y(t)=\psi(t)+A(t,T)\Th+B(t,T)\nu(t)+\int_t^T\Big[A(t,s)Y(s)+B(t,s)Z(s,t)\\
\ns\ds\qq\qq+D(t,s)Z(t,s)\Big]ds-\!\int_t^TZ(t,s)dW(s), \ t\in[0,T],\ \ \hb{a.e.} \\
\ns\ds \mu(t)=\Th\!-\!\int_t^T\nu(s)dW(s),\ \ \ \forall t\in[0,T],\\
\ns\ds \wt Y(t)\!=\!\wt \psi(t)\!+\!\int_t^T\!\big[\wt A(t,s)\wt Y(s)\!+\!
D(t,s)\wt Z(t,s)\big]ds\!-\!\int_t^T\!\wt Z(t,s)dW(s),\ \forall
t\in[0,T].
\ea\right.\ee

\mds

{\bf(H4)} Suppose $\Th\in L^2_{\cF_T}(\Omega;\dbR^m)$, $\wt\psi(\cd)\in C_{\cF_T}([0,T];L^2(\Omega;\dbR^m)),$ $\wt A:[0,T]^2\times\Omega\to\dbR^{m\times m}$
is measurable and bounded such that for almost all
$t\in[0,T]$, $s\1n\mapsto \wt A(t,s) $ is
$\dbF$-adapted on $[0,T]$. Moreover, for almost all $s\in[0,T]$,
the process $\wt A(\cd,s)$ is continuous.

\mds

As to the first equation of (\ref{BSVIE-adjoint-1}), under (H3)--(H4) there exists a unique pair of adapted M-solution $(Y(\cd),Z(\cd,\cd))\in \cH^2_1(0,T;\dbR^m)$ with (see \cite{Yong 2008})
$$
Y(t)=\dbE Y(t)+\int_0^tZ(t,s)dW(s),\qq \hb{a.e.}\ t\in[0,T],\ a.s.\
$$
As to the third one, by Lemma \ref{Pro-BSVIEs} there exists a unique solution $(\wt Y(\cd),\wt Z(\cd,\cd))\in \cH^2(0,T;\dbR^m)$. Now let us give the following kind of duality principle,

\mds

\bl\label{Lemma-dual}
\rm Suppose (H3)--(H4) hold true. Let $(\a,\b)\in C _\dbF([0,T];L^2(\Omega;\dbR^m))\times C([0,T];L^2_{\dbF}(0,T;\dbR^m))$ and $\xi$ be the solution to equation (\ref{SFVIE-1}). Then
\bel{duality-1}\q~~~\ba{ll}
\ns\ds\dbE\lan \xi(T),\Th\ran\!+\!\dbE\int_0^T\!\!\lan\psi(t),\xi(t)\ran dt
\!=\!\dbE\lan\alpha(T),\Th\ran\\
\ns\ds+\!\dbE\int_0^T\!\!\!\lan\beta(T,s),\nu(s)\ran ds\!+\!\dbE\int_0^T\!\!\!\lan Y(t),\a(t)\ran
dt\!+\!\dbE\int_0^T\!\!\!\int_0^T\!\!\lan Z(t,s),\b(t,s)\ran dsdt,
\ea\ee
and
\bel{duality-2}\ba{ll}
\ns\ds\q \dbE\lan \xi(0),\wt\psi(0)\ran\!+\!\dbE\int_0^T\!\!\lan
\dbE_s\xi(0),\wt A(0,s)\wt Y(s)\ran ds\\
\ns\ds=\!\lan \wt
Y(0),\alpha(0)\ran+\dbE\int_0^T\!\!\lan \wt Z(0,s),\beta(0,s)\ran ds.
\ea\ee
\el
\begin{proof} \rm Under (H3) it follows from Theorem \ref{SFVIE-regularity} that $\xi(0)$, $\xi(T)$
are well-defined. Given $\Th\in L^2_{\cF_T}(\Omega;\dbR^m)$, by the definition of $\xi(T)$,
$$\ba{ll}
\ns\ds \dbE\lan \xi(T),\Th\ran=\dbE\lan
\alpha(T),\Th\ran+\dbE\int_0^T\lan\beta(T,s),\nu(s)\ran
ds\\
\ns\ds\qq\qq\qq\ +\dbE\int_0^T\lan \xi(s), A(s,T)\Th+B(s,T)\nu(s)\ran ds,
\ea$$
On the other hand, by Theorem 4.1 in
\cite{Shi-Wang-Yong 2015}, one has,
$$\ba{ll}
\ns\ds\q \dbE\int_0^T\lan
\xi(s),\psi(s)+A(s,T)\Th+B(s,T)\nu(s)\ran ds\\
\ns\ds= \dbE\int_0^T\lan Y(s),\alpha(s)\ran ds+\dbE\int_0^T\int_0^T\lan
Z(t,s),\beta(t,s)\ran dsdt.
\ea
$$
As a result, by above two results we can obtain (\ref{duality-1}). Now let us treat the duality result
(\ref{duality-2}). Note that here $\wt Y(0)$ is well-defined.
On the other hand,
$$\ba{ll}
\ns\ds
\xi(0)=\alpha(0)+\int_0^T\Big[\beta(0,s)
+D(0,s)^{\top}\dbE_s \xi(0)\Big]dW(s).
\ea$$
Consequently,
$$\ba{ll}
\ns\ds \dbE\int_0^T\lan \wt Z(0,s),\beta(0,s)\ran ds=\dbE\lan \int_0^T \wt Z(0,s)dW(s),\int_0^T \beta(0,s)dW(s)\ran\\
\ns\ds=\dbE\lan \int_0^T \wt Z(0,s)dW(s),\xi(0)-\alpha(0)-\int_0^TD(0,s)^{\top}\dbE_s\xi(0)dW(s)\ran\\
\ns\ds=
\dbE\lan \int_0^T \wt Z(0,s)dW(s),\xi(0)\ran-\dbE\int_0^T\lan D(0,s)\wt Z(0,s),\xi(0)\ran ds.
\ea$$
Then, one can deduce that,
$$\ba{ll}
\ns\ds\q \dbE\lan \wt\psi(0),\xi(0)\ran\\
\ns\ds=\dbE\lan
\wt Y (0)-\int_0^T\big[\wt A(0,s)\wt Y(s)+D(0,s)\wt Z(0,s)\big]ds
+\int_0^T\wt Z(0,s)dW(s),\xi(0)\ran\\
\ns\ds =\dbE\lan \wt Y(0),\xi(0)\ran-\dbE\int_0^T\!\!\lan \wt
Y(s),\wt A(0,s)^{\top}\dbE_s\xi(0)\ran ds+\dbE\int_0^T\!\!\lan \wt
Z(0,s),\beta(0,s)\ran ds.
\ea$$
Note that $\dbE\xi(0)=\alpha(0)$, $\wt Y(0)$ is a constant, therefore one has (\ref{duality-2}).
\end{proof}

\br
\em
Let us make some points on the dualities established in Lemma \ref{Lemma-dual}. Firstly, (\ref{duality-1}) would degenerate into Theorem 4.1 of \cite{Shi-Wang-Yong 2015} when $\Th=0$, and Theorem 5.1 of \cite{Yong 2008} when $\beta(\cd,\cd)=0$. Secondly, if $\beta(\cd,\cd)=\Theta=0$, then (\ref{duality-1}) can also be obtained even when $\xi(\cd)\in L^2_{\cF_T}(0,T;\dbR^m)$, see Theorem 4.1 of \cite{Shi-Wang-Yong 2015}.
Thirdly, it seems that (\ref{duality-2}), which is
used to treat the term $\wt Y(0)$, $\wt Z(0,s)$, appears for the first time. In particular, when
$\beta(\cd,\cd)=0$, (\ref{duality-2}) plays the similar role as
Lemma 5.2 in \cite{Shi-Wang-Yong 2015}.
\er

\subsection{A pointwise procedure via set-valued analysis}
In this subsection, we introduce a technical lemma with the help of Proposition \ref{Pro-adjacent}, which is important in deriving the necessary conditions for optimal controls in the pointwise form. To fulfill the completeness requirement there, we need some preparations.

Firstly, let us recall some notions about the set-valued stochastic processes, we refer the reader to \cite{Kisielewicz13} for more details. We call a measurable set-valued map $\zeta:(\Omega,\cF)\to 2^{\dbR^l}$ a set-valued random variable, and a family of set-valued random variables $\Gamma(t,\cdot): \Omega\to 2^{\dbR^l}$, $t\in[0,T]$ a set-valued stochastic process. $\Gamma$ is called to be measurable, if it is $\mb([0,T])\otimes\cF$-measurable. Furthermore, if $\Gamma(t,\cdot)$ is $\cF_{t}$-measurable for any $t\in [0,T]$, $\Gamma$ is named as $\dbF$-adapted.
Define:
$$ \cG:=\big\{A\in\cB([0,T])\times\cF,\big|\ A_t\in\cF_t,\
\forall t\in[0,T]\big\},\qq A_t:=\{\o\in\Omega\ \big|\ (t,\o)\in A\}.
$$
It is easy to see that $\cG\subset \cB([0,T])\times\cF$ is a
$\sigma$-algebra and $([0,T]\times\Omega,\cG,\l\times\dbP)$ is a measure
space. Moreover, the following result holds true, see p. 96 of \cite{Kisielewicz13}.

\mds

\begin{lemma}\label{lemma adapted sigma field}
\rm Let $\Gamma:([0,T]\times\Omega,\cB([0,T])\otimes\cF)\to 2^{\mrm}$. $\Gamma$ is measurable and $\dbF$-adapted if and only if $\Gamma$ is $\cG$-measurable.
\end{lemma}

\mds

By Lemma \ref{lemma adapted sigma field}, it is easy to see that
$$
L^2_{\dbF}(0,T;\dbR^l)=\Big\{y:[0,T]\times\Omega\rightarrow \dbR^l\ \big|\ y(\cd)\
\hb{is
$\cG$-measurable,}\ \dbE\int_0^T|y(s)|^2ds<\infty\Big\},$$
$$\ba{ll}
\ns\ds\cU_{ad}\!=\!\Big\{u:[0,T]\!\times\!\Omega\to \dbR^l\ \big|\ u(\cd)\
\hb{is $\cG$-measurable,}\\
\ns\ds\qq\qq\qq\qq\qq\qq u(t)\!\in\! U,\  \hb{a.e.} \ (t,\omega)\in[0,T]\!\times\!\Omega, \
\dbE\int_0^T\!|u(s)|^2ds<\infty\Big\}.
\ea$$
Note that $([0,T]\times\Omega;\cG,\l\times\dbP)$ may not be complete, hence in the following we need to construct a
complete version. The $completion$ of
$\cG$ under $\lambda\times\dbP$, i.e. $\cG^*$ is the collection of
subset $A$ of $[0,T]\times\Omega$ for which there exist $E,\ F\in\cG$
such that
$$E\subset A\subset F,\ \hb{and} \ [\lambda\times\dbP](F-E)=0.$$
In
this case one can define a function $\mu^*$ on $\cG^*$ as
$\mu^*(A)=[\lambda\times\dbP](E)$. According to Proposition 1.5.1 in
\cite{Cohn13},  $([0,T]\times\Omega;\cG^*,\mu^*)$ is a completion of $([0,T]\times\Omega;\cG,\lambda\times\dbP)$. Moreover, it also implies the following result which is useful later.

\mds

\begin{lemma}\label{complete-existence-a.s.}
 \rm Let $(\Xi, \ms, \mu)$ be a $\sigma$-finite measure space with its completion $(\Xi,
\ms^*, \mu^*)$, $f$ be a $\ms^*$-measurable function from $\Xi$ to
$\dbR^l$. Then there exists a $\ms$-measurable function $g$ such
that $\mu^*(g(\xi)\neq f(\xi))=0$.
\end{lemma}

\mds

\rm

Define
\bel{cG=cF-1-*}\qq\ba{ll}
\ns\ds \cL^{2}_{\dbF}(0,T;\dbR^l)\!:=\!\Big\{y:[0,T]\!\times\!\Omega\rightarrow \dbR^l\ \big|\
y(\cd)\ \hb{is $\cG^*$-measurable,}\\
\ns\ds\qq\qq\qq\qq\qq\qq\qq\qq\qq \int_{[0,T]\!\times\!\Omega}|y(s,\o)|^2d\mu^*(s,\o)\!<\!\infty\Big\},\\
\ns\ds \cU^*_{ad}\!:=\!\Big\{u:[0,T]\!\times\!\Omega\rightarrow \dbR^l\ \big|\ u(\cd)\
\hb{is $\cG^*$-measurable,}\ u(t)\!\in\! U,\\
\ns\ds\qq\qq\qq\qq\qq\qq\qq \mu^*\!-\!\hb{a.e.}, \
\int_{[0,T]\!\times\!\Omega}|u(s,\o)|^2d\mu^*(s,\o)\!<\!\infty\Big\}.
\ea\ee
Clearly,
$\cU_{ad}\subset\cU^*_{ad}$ and $L^{2}_{\dbF}(0,T;\dbR^l)\subset \cL^{2}_{\dbF}(0,T;\dbR^l).$
In particular, if we suppose that,
$$\Xi=[0,T]\times\Omega,\ \ \ms=\cG^*,\ \ \mu=\mu^*,\ \ X=\dbR^l,$$
by
Proposition \ref{Pro-adjacent} one has,

\medskip
\bl \label{Lem-adjacent}  \rm Let $U$ be a closed subset of $\dbR^l$.
Then, for any $u(\cdot)\in \mmu^*$, $\T_{U}(u(\cdot)):[0,T]\times\Omega\rightsquigarrow \dbR^l$ is
$\cG^*$-measurable, and $\mt_u^* \subset \T_{\mmu^*}(u(\cdot))$,
where
\bel{cT-bar-u*}\ba{ll}
\ns\ds \mt_u^*:=\big\{v(\cdot)\in \cL^{2}_{\dbF}(0,T;\dbR^l)\ \big|\
v(t)\in \T_{U}(u(t)),\ u^*-\hb{a.e.}\big\}.
\ea\ee
\el

\mds

Now, we give the main result of this subsection.

\mds

\bl\label{Lemma-integral-pointwise}
\rm Suppose $\bar u(\cd)$ is an optimal control,
$F:\Omega\times[0,T]\rightarrow\dbR^l$ is a
$\cB([0,T])\times\cF$-measurable and $\dbF$-adapted process such
that,
\bel{integral-v}\ba{ll}
\ns\ds \me\int_{0}^{T}\lan F(t),v(t)\ran dt\geq 0,\qq \forall
v(\cdot)\in \T_{\mmu}(\bu(\cdot)).
\ea\ee
Then we have,
\bel{pointwise-v}\ba{ll}
\ns\ds \lan F(t),v \ran \geq0,\ \ \ \forall v\in T^{b}_{U}(\bar
u(t)),\ \ \t\in[0,T].\ \ \  [\l\times\dbP]-\hb{a.e.}
\ea\ee
\el
\begin{proof} \rm We would like to slip the proof into several parts.

\it Step 1: \rm In this step, we prove that
\bel{integral-v*}\ba{ll}
\ns\ds \int_{[0,T]\times\Omega}\lan F(t,\omega),v^*(t,\omega)\ran
d\mu^*(t,\omega)\geq 0,\qq \forall v^*(\cdot)\in
\T_{\mmu^*}(\bu(\cdot)),
\ea\ee
with $\cU_{ad}^*$ defined by (\ref{cG=cF-1-*}).

For any $v^*(\cd)\in \T_{\mmu^*}(\bu(\cdot))$, we know that
$v^*(\cd)\in \cL^{2}_{\dbF}(0,T;\dbR^l)$. By Lemma
\ref{complete-existence-a.s.}, there exists a $\cG$-measurable
function $v(\cd)$ on $[0,T]\times\Omega$ such that
\bel{v*-v}\q\ba{ll}
\ns\ds v^*(s,\o)=v(s,\o),\ \mu^*\!-\!a.e. \ \Rightarrow
\int_{[0,T]\times\Omega}\!|v^*(s,\o)-v(s,\o)|^2d\mu^*(s,\o)=0.\ \
\ea\ee
As a result, one has
$$
\dbE\int_0^T|v(s,\o)|^2ds=\int_{[0,T]\times\Omega}|v(s,\o)|^2d\mu^*(s,\o)<\infty,
$$
which implies that $ v(\cd)\in L^2_{\dbF}(0,T;\dbR^l)$. On the other
hand, by $v^*(\cd)\in \T_{\mmu^*}(\bu(\cdot))$ and Lemma
\ref{equivdef for adjacent cone}, for any $h_n\to 0^+$, there
exist $v^*_n(\cd)\in \cL^{2}_{\dbF}(0,T;\dbR^l)$ such that $\bar u(\cd)+h_nv^*_n(\cd)\in \cU_{ad}^*$, $n\in \mn$, and
\bel{v*-n-v*}\ba{ll}
\ns\ds
\int_{[0,T]\times\Omega}|v^*_n(s)-v^*(s)|^2d\mu^*(s,\o)\rightarrow 0,\ \
\ n\rightarrow\infty.
\ea\ee
Similar to the above, for any $n$,  there exists a
process $v_n(\cd)\in L^2_{\dbF}(0,T;\dbR^l)$ such that,
\bel{v*-n-v-n}\q~~\ba{ll}
\ns\ds v^*_n(s,\o)\!=\!v_n(s,\o),\ \mu^*\!-\!\hb{a.e.,} \   \hb{and}\
 \int_{[0,T]\times\Omega}\!|v^*_n(s,\o)-v_n(s,\o)|^2d\mu^*(s,\o)\!=\!0.
\ea\ee
Combining
(\ref{v*-v}), (\ref{v*-n-v*}) with (\ref{v*-n-v-n}) one has,
$$ \dbE\int_0^T\!|v_n(s)\!-\!v(s)|^2ds\!=\!
\int_{[0,T]\times\Omega}|v_n(s,\o)\!-\!v(s,\o)|^2d\mu^*(s,\o)\!\to\! 0,\
\ n\!\to\!\infty,
$$
and $\bar u(\cd)+h_nv_n(\cd)\in \cU_{ad}$.
Then, by Lemma \ref{equivdef for adjacent cone},  $v(\cd)\in
T^b_{\cU_{ad}}(\bar u(\cd))$ and
\begin{eqnarray*}
&&  \int_{[0,T]\times\Omega}\!\lan F(t,\o),v^*(t,\o)\ran d\mu^*(t,\o)\\
&=&
\int_{[0,T]\times\Omega}\!\lan F(t,\o),v(t,\o)\ran
d\mu^*(t,\o)\!=\!\dbE\int_0^T\!\lan F(t),v(t)\ran dt\ge0.
\end{eqnarray*}
This proves (\ref{integral-v*}).

\it Step 2: \rm In this step, we prove the set
\bel{cA-bar-u}\ba{ll}
\ns\ds \cA_{\bar u}:=\big\{(t,\omega)\in[0,T]\times\Omega\ \big|\
\lan F(t),v\ran \ge 0, \ \ \forall\ v\in \T_{U}(\bu(t))\big\}
\ea\ee
is $\cG^*$-measurable.
Let us first look at its complement, i.e.
\bel{cA-bar-u-c}\ba{ll}
\ns\ds \cA^c_{\bar u}:=\big\{(t,\omega)\in[0,T]\times\Omega\big|\
\exists \ v\in \T_{U}(\bu(t))\  \hb{s.t.} \  \inner{F(t)}{v}<0\big\}.
\ea\ee
By Lemma \ref{Lem-adjacent} above the set-valued map $\T_{U}(\bar
u(\cd)):[0,T]\times\Omega\rightsquigarrow \dbR^l$ is $\cG^*$-measurable.
Hence according to Proposition \ref{Pro-measurable},
$$
\big\{(t,\o,v)\in[0,T]\times\Omega\times\dbR^l \big|\ v\in\T_U(\bar
u(t,\o))\big\}\in\cG^*\otimes\cB(\dbR^l).
$$
By the assumption on $F(\cd)$, we have,
\bel{Graph-1}\ba{ll}
\ns\ds
 \big\{(t,\o,v)\in[0,T]\times\Omega\times\dbR^l \big|\
v\in\T_U(\bar u(t,\o)),\ \inner{F(t)}{v}<0
\big\}\in\cG^*\otimes\cB(\dbR^l).
\ea\ee
Now let us define a set-valued map
$\L(\cd,\cd):[0,T]\times\Omega\rightsquigarrow \dbR^l$ as,
$$
\L(t,\o):=\big\{v\in \dbR^l \big|\ v\in\T_U(\bar u(t,\o)),\
\inner{F(t)}{v}<0 \big\},\qq (t,\o)\in[0,T]\times\Omega.
$$
By Proposition \ref{Pro-measurable}, it follows from
(\ref{Graph-1}) that map $\L$ is $\cG^*$-measurable. Then $\cA^c_{\bar u}$, the
domain of map $ \L $, is measurable. Consequently,
$\cA_{\bar u}$ is $\cG^*$-measurable.

\it Step 3: \rm In this step we  prove that
$\mu^*(\cA^c_{\bar u})=0$.

Let $k,\
r=1,2,\cdots$,  define
\begin{eqnarray*}
\cA_{\bu}^{k,r}:=\Big\{(t,\omega)\in [0,T]\times \Omega&\Big|&
\exists\ v\in \T_{U}(\bu(t))\cap \bar B(0,r),\ s.t.\ \
\inner{F(t)}{v } \le -\frac{1}{k} \Big\}.
\end{eqnarray*}
It is clear that
$$ \cA^c_{\bar u}=\bigcup_{k\geq1}\bigcup_{r\geq1}\cA_{\bu }^{k,r}.
$$
Like above $\cA_{\bar u}$ one can prove that $\cA_{\bu }^{k,r}$ is
$\cG^*$-measurable which implies that $(\cA_{\bu }^{k,r},\cG^*)$ is
a measurable space. If for any $k,\ r\geq 1$, we can prove $\cA_{\bu
}^{k,r}$ has zero measure, then the proof can be finished. We will
prove this by contradiction. Suppose that there exist $k$ and $r$ such that $\mu^*(\cA_{\bu }^{k,r})>0$.  Define a set-valued map $\G:\cA_{\bu
}^{k,r}\rightsquigarrow\dbR^l$ by
$$\Gamma^{k,r}(t,\omega):=\Big\{v\in \T_{U}(\bu(t))\cap \bar
B(0,r)\ \Big|\ \inner{F( t)}{v } \le -\frac{1}{k} \Big\}.$$
Obviously, $\G^{k,r}(t,\o)$ is closed-valued. Similar as (\ref{Graph-1}), the set
\bel{}\ba{ll}
\ns\ds
 \big\{(t,\o,v)\in[0,T]\times\Omega\times\dbR^l \big|\
v\in\T_U(\bar u(t,\o))\cap \bar B(0,r),\ \inner{F(t)}{v}\leq -\frac
1 k \big\},
\ea\ee
is $\cG^*\otimes\cB(\dbR^l)$-measurable, from which, as well as Proposition \ref{Pro-measurable} one can
obtain that $\Gamma^{k,r}_{\bar u}$ is a $\cG^*$-measurable
set-valued map with $\hb{Dom}(\Gamma^{k,r})=\cA_{\bu }^{k,r}$. Then by Proposition \ref{Pro-mea-sel} there exists
a $\cG^*$-measurable selection
$v^{k,r}(\cdot)$ on $\hb{Dom}(\Gamma^{k,r})$,
 i.e.,
$$v^{k,r}(t,\omega)\in \Gamma^{k,r}(t,\omega)\subset\big[ \T_{U}(\bu(t))\cap
\bar B(0,r)\big], \ \forall\ (t,\omega)\in
\hb{Dom}(\Gamma^{k,r})=\cA^{k,r}_{\bar u }.$$
Define $\wt v^{k,r}(\cd):= v^{k,r}(\cd
)I_{\mt^{k,r}_{\bar u }}(\cd )$, then $\wt
v^{k,r}(\cd)\in\cA^*_{\bar u}$, where $\mt^*_{\bar u}$ is defined by
(\ref{cT-bar-u*}), and
\bel{negative-1}\ba{ll}
\ns\ds\mu^*\big\{( t,\omega)\in[0,T]\big|\ \inner{F( t)}{\wt
v^{k,r}(t)} \le-\frac{1}{k}\big\}\ge \mu^*(\cA^{k,r}_{\bar u})>0.
\ea\ee
Consequently,
\begin{equation}\label{add eq4.27}
\int_{[0,T]}\lan F(t,\o),\wt v^{k,r}(t,\o)\ran
d\mu^*(t,\o)\le -\frac{1}{k}\mu^*(\cA^{k,r}_{\bar u})<0.
\end{equation}

On the other hand, by Lemma \ref{Lem-adjacent} one has $v^{k,r}(\cdot)\in\mt^*_{\bar u}\subset\T_{\cU_{ad}^*}(\bar u(\cd))$. It then follows from
(\ref{integral-v*}) that
$$\int_{[0,T]}\lan F(t,\o),\wt v^{k,r}(t,\o)\ran
d\mu^*(t,\o)\ge0,$$
which contradicts to (\ref{add eq4.27}). Therefore,
$\cA_{\bu}^{k,r}$ has zero measure.

\it Step 4: \rm In this step we would like to prove that there
exists $\cG$-measurable set $\cD_{\bar u}\subset\cA_{\bar u}$
satisfying $[\lambda\times\dbP](\cD_{\bar u})=T$ which naturally
implies the conclusion.

Actually, for above $\cA^c_{\bar u}\in\cG^*$, by the definition of
$\cG^*$, there exists a $\cG$-measurable set $\cE_{\bar u}$
satisfying $\cA^c_{\bar u}\subset \cE_{\bar u}$ and
$\mu^*(\cA^c_{\bar u})=\mu(\cE_{\bar u})=0$. For this
$\cE_{\bar u}$, by \it Step 3 \rm one immediately has $\cD_{\bar
u}:= \cE_{\bar u}^c \subset\cA_{\bar u}$ and $[\lambda\times\dbP](\cD_{\bar u})=T$.
Then the proof is finished.
\end{proof}

\subsection{Proofs of Theorem \ref{Theorem-main-result}}

In this part, based on the preparations in the last two subsections, we
are about to give the proof of the main result in Section 3.

Firstly we need to introduce the variational equations.
Let $\bar{u}(\cdot)$ be an optimal control and
$v(\cdot)\in \T_{\mmu}(\bu(\cdot))$. By Remark \ref{remark2-1}, for
any $\eps>0$ there exists a $v_{\eps}(\cdot)\in
L_{\dbF}^{2}(0,T;\dbR^l)$ such that $u^{\e}(\cd):=\bu(\cdot)+\eps
v_{\eps}(\cdot)\in \mmu$ and
$$\big\|v-v_{\e}\big\|^2_{L^2_{\dbF}(0,T;\dbR^l)}:=\me\int^{T}_{0}
|v(t)-v_{\eps}(t)|^2dt\to 0,\ \eps\to 0^+.$$
As a result, there exists a constant $\delta<1$ such that
$$\sup_{\e\in(0,\delta]}\big\|v_{\e}\big\|_{L^2_{\dbF}(0,T;\dbR^l)}
:=\sup_{\e\in(0,\delta]}\me\int^{T}_{0}|v_{\eps}(t)|^2dt<\infty.$$
Suppose $(X^{\e}(\cd),Y^{\e}(\cd),Z^{\e}(\cd,\cd))\ (\bar
X(\cd),\bar Y(\cd),\bar Z(\cd,\cd))$ are the state processes
associated with $u^{\e}(\cd),\ \bar u(\cd)$ respectively.

Given FBSVIE (\ref{FBSVIE-3-1}), we introduce the following first order variational
system,
\bel{variational-1}\left\{\2n\ba{ll}
\ns\ds X_1(t)\!=\!\!\int_0^t\!\!\big[b_x(t,s)X_1(s)\!+\! b_u(t,s)v(s)\big]ds
\!+\!\!\int_0^t\!\!\big[\si_x(t,s)X_1(s)\!+\!\si_u(t,s)v(s)\big]dW(s),\\
\ns\ds Y_1(t)\!=\!\psi_x(t)X_1(T)\!+\!\int_t^T\!\!\Big\{g_x(t,s)X_1(s)\!+\!g_y(t,s)Y_1(s),\\
\ns\ds\qq\qq+g_z(t,s)Z_1(t,s)\1n+\1n g_u(t,s)v(s)
\Big\}ds\1n-\2n\int_t^T\2n Z_1(t,s)dW(s),\qq  t\in[0,T],
\ea\right.\ee
where e.g. $b_x(t,s):= b_x(t,s,\bar X(s),\bar u(s))$.
Under (H1), we know that (\ref{variational-1}) admits a unique
adapted solution $(X_1(\cd),Y_1(\cd),Z_1(\cd\,,\cd))$ in the sense
of Definition \ref{Def-M-SVIEs} and Definition \ref{Def-M-BSVIEs}.
For $t,s\in [0,T]$, define
\bel{5.5}\q
X_1^\e(t)\!:=\!\frac{X^\e(t)\!-\!\bar X(t)}{\e},\
Y^\e_1(t)\!:=\!\frac{Y^\e(t)\!-\!\bar Y(t)}{\e},\
Z_1^\e(t,s)\!:=\!\frac{Z^\e(t,s)\!-\!\bar Z(t,s)}{\e}.\ee
\bl\label{convergence-1}
 \rm Suppose (H1) hold true, $(X_1^{\e}(\cd),Y_1^{\e}(\cd),Z_1^{\e}(\cd,\cd))$ is the state
processes associated with $u^{\e}(\cd)$, $(X_1 (\cd),Y_1 (\cd),Z_1
(\cd,\cd))$ is the unique C-adapted solution of FBSVIE
(\ref{variational-1}) in the sense of Definition \ref{Def-M-SVIEs} and Definition \ref{Def-M-BSVIEs}. Then we have,

\bel{convergence-result-1}\q \ba{ll}
\ns\ds\lim_{\e\to0}\dbE|X_1^\e(t)-X_1(t)|^2\!=\!0,\    \forall t\in[0,T],\ \ \lim_{\e\to0}\dbE\int_0^T\!|X_1^\e(t)-X_1(t)|^2dt\!=\!0,
\ea \ee
and
\bel{convergence-result-2}~~ \ba{ll}
\ns\ds\lim_{\e\to0}\(\dbE|Y^\e_1(t)-Y_1(t)|^2+\dbE\int_t^T|Z_1^\e(t,s)-Z_1(t,s)|^2ds\)=0,\ \ \forall
t\in[0,T].
 \ea\ee
\el
\begin{proof} \rm By the
standard estimates for SVIEs and BSVIEs in Lemma
\ref{Pro-BSVIEs}, one has,
\bel{estimate-1}\qq\left\{\ba{ll}
\ns\ds \sup\limits_{t\in[0,T]}\dbE|\bar
X(t)|^2\!+\!\sup\limits_{t\in[0,T]}\dbE|X^{\e}(t)|^2\!+\!\sup\limits_{t\in[0,T]}\dbE|\bar
Y(t)|^2\!+\!\sup\limits_{t\in[0,T]}\dbE|Y^{\e}(t)|^2
<\infty,\\
\ns\ds \sup\limits_{t\in[0,T]}\dbE\int_t^T|\bar
Z(t,s)|^2ds+\sup\limits_{t\in[0,T]}\dbE\int_t^T|Z^{\e}(t,s)|^2ds<\infty,
 \q \forall\e\in(0,\delta],\\
\ns\ds \sup_{t\in[0,T]}\dbE|X^\e(t)\1n-\1n \bar X(t)|^2\1n  \le
C\e^2\dbE\int_0^T\2n| v_{\e}(s)|^2ds,\\
\ns\ds \sup_{t\in[0,T]}\dbE|Y^\e(t)\1n - \1n\bar Y(t)|^2
\1n\!+\!\2n\sup_{t\in[0,T]}\dbE\2n\int_t^T \2n|Z^\e(t,s)\1n - \1n\bar
Z(t,s)|^2ds\!\le\! C\e^2\dbE\int_0^T\2n| v_{\e}(s)|^2ds.
\ea\right.\ee

Similarly,
$$\ba{ll}
\ns\ds
\sup_{t\in[0,T]}\dbE|X_1(t)|^2\!+\!\sup_{t\in[0,T]}\dbE|Y_1(t)|^2\!
+\!\sup_{t\in[0,T]}\dbE\int_t^T|Z_1(t,s)|^2ds\le
C\dbE\int_0^T|v(t)|^2dt.
\ea$$

It is a simple fact that $ X^\e_1(\cd)-X_1(\cd)$ is the solution to the following SVIE:
$$\ba{ll}
\ns\ds
\q X^\e_1(t)-X_1(t)\\
\ns\ds=\int_0^t\Big\{b_x^\e(t,s)X^\e_1(s)-b_x(t,s)X_1(s)+\big[
b_u^\e(t,s)v_{\e}(s)-b_u(t,s)v(s)\big]\Big\}ds\\
\ns\ds\q +\!\int_0^t\Big\{\si_x^\e(t,s)X^\e_1(s)\!-\!\si_x(t,s)X_1(s)
\!+\!\big[\si_u^\e(t,s)v_{\e}(s) \!-\!\si_u(t,s)v(s)\big]\Big\}dW(s),\ea$$
where,
$$b_x^\e(t,s):=\int_0^1b_x\big(t,s,\bar X(s)+\th[X^\e(s)-\bar X(s)],\bar
u(s)+\th\e v_{\e}(s)\big)d\th,$$
and $b^\e_u(t,s)$, $\si^\e_x(t,s)$, $\si^\e_u(t,s)$ are defined in a
similar manner. As a result,
it follows from dominated convergence theorem and Lemma
\ref{Pro-BSVIEs} that (\ref{convergence-result-1}) holds true. To obtain similar convergence result for the backward equation, let
us look at,
$$\ba{ll}
\ns\ds Y^\e_1(t)-Y_1(t)=\psi^\e_x(t)X_1^\e(T)-\psi_x(t)X_1(T)
+\int_t^T\Big\{g_x^\e(t,s)X_1^\e(s)-g_x(t,s)X_1(s)\\
\ns\ds\q+g_y^\e(t,s)Y_1^\e(s)-g_y(t,s)Y_1(s)+g_z^\e(t,s)Z_1^\e(t,s)-g_z(t,s)Z_1(t,s)\\
\ns\ds\q+\big[g_u^\e(t,s)v_{\e}(s)-g_u(t,s)v(s)\big]\Big\}ds-\int_t^T\(Z_1^\e(t,s)-Z_1(t,s)\)dW(s),\ea$$
where for example,
$$\psi_{x}^\e(t)=\int_0^1\psi_{x}\big(t,\bar X(T)+\th[X^\e(T)-\bar X(T)]\big)
d\th,$$
Then by dominated convergence theorem and the estimates in Lemma
\ref{Pro-BSVIEs} we have (\ref{convergence-result-2}). This completes the proof of Lemma \ref{convergence-1}.
\end{proof}

\mds

Now, we give the proof of Theorem  \ref{Theorem-main-result}.

\mds

\it\bf Proof of Theorem \ref{Theorem-main-result}. \rm By the optimality of $(\bar
X(\cd),\bar Y(\cd),\bar Z(\cd\,,\cd),\bar u(\cd))$,
\bel{J-1-1}\ba{ll}
\ns\ds0\le\frac{J(u^\e(\cd))-J(\bar u(\cd))}{\e}\\
\ns\ds\ \ \ =\dbE\big[\lan h^\e_x,X^\e_1(T)\ran+\lan h^\e_y, Y^\e_1(0)\ran\big]+\dbE\! \int_0^T\!\!
\big[\lan f_x^\e(s),X^\e_1(s)\ran\!\\
 \ns\ds\qq\qq\qq+\!\lan f_y^\e( s),Y^\e_1(s)\ran \!+\!\lan f_z^\e(
s),Z^\e_1(0,s)\ran\!+\!\lan f_u^\e(s),v_{\e}(s)\ran\big]ds,\ea\ee
where for example,
$$\ba{ll}
\ns\ds h^\e_x:=\int_0^1h_x(\bar X(T)+\th[X^\e(T)-\bar X(T)],\bar
Y(0)+\th[Y^\e(0)-\bar Y(0)])d\th,\\
\ns\ds f^\e_x( s)=\int_0^1f_x( s,\wt X(s;\th),\wt Y(s;\th),\wt Z(0,s;\th),
\wt u(s;\th))d\th,\\
\ns\ds \wt \varphi(s;\th)=\bar \varphi(s)+\th[\varphi^\e(s)-\bar \varphi(s)],\ \ \varphi=X,\ Y,\ Z,\ u.
\ea$$
By (\ref{convergence-result-1})--(\ref{estimate-1}), passing
to the limit in (\ref{J-1-1}), we obtain that
\bel{J-1-2}\ba{ll}
\ns\ds0\le \dbE\big[\lan h_x,X_1(T)\ran+\lan h_y, Y_1(0)\ran\big]+\dbE \int_0^T
\big[\lan f_x(s),X_1(s)\ran\\
 \ns\ds\qq\qq\q +\lan f_y( s),Y_1(s)\ran +\lan f_z(
s),Z_1(0,s)\ran+\lan f_u(s),v(s)\ran\big]ds.
\ea\ee
Given inequality (\ref{J-1-2}), next we are about to obtain the
following necessary condition in the integral form,
\bel{cH-1-inequality}\ba{ll}
\ns\ds \dbE\int_0^T\lan \mH_{u}(s;\bar X(\cd),\bar Y(\cd),\bar
Z(\cd,\cd),\bar u(\cd)),v(s)\ran ds\geq0, \q \forall\ v(\cd)\in T^b_{\cU_{ad}}(\bar u(\cd)),
\ea\ee
where $\mH_{u}(s;\bar X(\cd),\bar Y(\cd),\bar
Z(\cd,\cd),\bar u(\cd))$ is defined by (\ref{cH-1}). Then the pointwise necessary condition (\ref{maximum-1}) follows
from Lemma \ref{Lemma-integral-pointwise}. To obtain
(\ref{cH-1-inequality}), firstly let us deal with $Y_1(0)$, $Z_1(0,\cd)$ in
(\ref{J-1-2}). Note that
$$\ba{ll}
\ns\ds
Y_1(0)=\wt\Psi+\int_0^T\big[g_y(0,s)Y_1(s)+g_z(0,s)Z_1(0,s)\big]ds
-\int_0^TZ_1(0,s)dW(s),
\ea$$
where
\bel{wt-psi-1}\ba{ll}
\ns\ds \wt\Psi:=
\psi_x(0)X_1(T)+\int_0^T\big[g_x(0,s)X_1(s)+g_u(0,s)v(s)\big]ds,
\ea\ee
and $v(\cdot)\in \T_{\mmu}(\bar u(\cd))$.
According to Theorem \ref{SFVIE-regularity}, there exists a unique
$\cF_T$-measurable random variable $\l(0)$ satisfying the following
equality,
$$\l(0)=\dbE h_y(\bar X(T),\bar
Y(0))^{\top}+\int_0^T\big[f_z(s)+g_z(0,s)\dbE_s\l(0)\big]dW(s).
$$
Denote $\alpha(0):=\dbE h_y(\bar X(T),\bar Y(0))^{\top},$ $\beta(0,\cd):=f_z(\cd)^{\top}$,
$\wt\psi(0):=\wt\Psi$, $\wt A(0,s):=g_y(0,s)$. By the duality result in
(\ref{duality-2}), we have
\bel{add 4.39}\ba{ll}
\ns\ds
\q\dbE\lan h_y,Y_1(0)\ran+\dbE\int_0^T\lan Z_1(0,s),f_z(s)\ran ds\\
\ns\ds=\dbE\lan\l(0),\wt\Psi\ran+\dbE\int_0^T\lan
\dbE_s\l(0),g_y(0,s)Y_1(s)\ran ds\\
\ns\ds=\!\dbE\lan\l(0),\psi_x(0)X_1(T)\ran\!+\dbE\int_0^T\!\!\lan
\dbE_s\l(0),g_x(0,s)X_1(s)\!\\
\ns\ds\qq\qq\qq\qq\qq\qq\qq+\!g_y(0,s)Y_1(s)\!+\!g_u(0,s)v(s)\ran ds.
\ea\ee
Substituting (\ref{add 4.39})  into (\ref{J-1-2}) we then obtain that
\bel{J-1-variational}\q\ba{ll}
\ns\ds 0\!\leq\!\dbE\Big\{\lan h_x\!+\!\psi_x(0)^{\top}\l(0),X_1(T)\ran
\!+\!\int_0^T\!\!\[\lan g_x(0,s)^{\top}\dbE_s\l(0)\!+\!f_x(s),X_1(s)\ran
\\
\ns\ds\q\ +\!\lan g_u(0,s)^{\top}\dbE_s\l(0)\!+\!f_u(s),v(s) \ran\]ds\!+\!\int_0^T\!\!\2n\lan
g_y(0,s)^{\top} \dbE_s\l(0)\!+\!f_y(s),Y_1(s)\ran ds\Big\}.\ea\ee
Now let us turn to deal with $Y_1(\cd)$ by means of Lemma
\ref{Lemma-dual}. To this end, we choose,
$$\left\{\ba{ll}
\ns\ds A(t,s):= g_y(t,s),\ \ B(t,s):= g_z(t,s),\ \ C(t,s):= 0,\ \ \Th:= 0,\ \
 \beta:= 0,\\
\ns\ds\a(t):= g_y(0,t)^{\top}\dbE_t\l(0)+f_y(t),\\
\ns\ds
\psi(t):= \psi_x(t)X_1(T)+\int_t^T\big[ g_x(t,s)X_1(s)
+g_u(t,s) v(s)\big]ds.
\ea\right.
$$
Then by the duality result (\ref{duality-1}),
\begin{eqnarray*}
&& \dbE \int_0^T\lan g_y(0,t)^{\top}\dbE_t\l(0)+f_{y}(t),Y_1(t)\ran dt\\
&=&\dbE\int_0^T\lan\psi_x(t)X_1(T)+\int_t^Tg_x(t,s)X_1(s)ds
+\int_t^Tg_u(t,s) v(s)ds,\xi(t)\ran dt.
\end{eqnarray*}
Note that here $\xi(\cd)\in L^2_{\cF_T}(0,T;\dbR^m)$. As a result, (\ref{J-1-variational}) can be rewritten as,
\bel{J-1-variational-2}\ba{ll}
\ns\ds 0\leq\dbE\Big\{\lan
h_x+\psi_x(0)^{\top}\l(0)+\int_0^T\psi_x(t)^{\top}\xi(t)dt,X_1(T)\ran\\
\ns\ds\qq\qq+\int_0^T\[\lan
g_x(0,s)^{\top}\dbE_s\l(0)+f_x(s)+\int_0^sg_x(t,s)^{\top}\xi(t) dt,X_1(s)\ran\\
\ns\ds\qq\qq\qq+\lan
g_u(0,s)^{\top}\dbE_s\l(0)+\int_0^sg_u(t,s)^{\top}\xi(t)dt+f_u(s),v(s)
 \ran\]ds\Big\}.
 \ea\ee
At last let us deal with the term for $X_1(T)$, $X_1(\cd)$ by means
of Lemma \ref{Lemma-dual} again. To this end, let us denote by
$$\ba{ll}
\ns\ds \Th:=
h_x+\psi_x(0)^T\l(0)+\int_0^T\psi_x(r)^{\top}\xi(r)dr,\qq
\beta(\cd,\cd):=0,\\
\ns\ds \psi(s):= g_x(0,s)^{\top}\dbE_s\l(0)+\int_0^sg_x(t,s)^{\top}\xi(t)
dt+f_x(s)^{\top},\\
\ns\ds\alpha(t):=
\int_0^tb_u(t,s)v(s)ds+\int_0^t\sigma_u(t,s)v(s)dW(s).
\ea$$
Then by duality (\ref{duality-1}) we have,
$$\ba{ll}
\ns\ds\q \dbE\lan X_1(T),\Th\ran+\dbE\int_0^T\lan \psi(t),X_1(t)\ran
dt\\
\ns\ds=\dbE\int_0^T\lan
b_u(T,s)v(s),\Th\ran ds+\dbE\int_0^T\lan \si_u(T,s)v(s),\nu(s)\ran ds\\
\ns\ds\q+\dbE\int_0^T\lan p(t),\int_0^tb_u(t,s)v(s)ds\ran
dt+\dbE\int_0^T\int_0^t\lan q(t,s),\si_u(t,s)v(s) \ran ds dt,
\ea$$
where $(p(\cd),q(\cd,\cd))$ is the solution to the last equation of the adjoint system (\ref{adjoint-1}).
Then, one can rewrite (\ref{J-1-variational-2}) as
(\ref{cH-1-inequality}), and the integral type condition (\ref{cH-1-inequality}) holds. This finishes the proof of Theorem \ref{Theorem-main-result}.

\section{Concluding remark}

This paper is devoted to optimal control problems of forward-backward stochastic Volterra integral equations. The control region is supposed to be closed, but not necessary to be convex. A new variational approach is introduced, which enables us to drop the second-order adjoint equations and weaken the regularity assumption on the involved coefficients. Note that these ideas are even new under the special stochastic differential equations, forward-backward stochastic differential equations and forward stochastic Volterra integral equations frameworks.
However, just as the above Remark \ref{strucural remark} shows, the variational approach requires some relatively strong structural assumptions on the control regions. Therefore, it still remains its importance to establish a general stochastic maximum principle for optimal control problems of forward-backward stochastic Volterra integral equations, especially when the control regions do not satisfies the structural assumption (\ref{strucural assum}). We will discuss this topic in our forthcoming paper.

\end{document}